\documentclass[12pt]{article}
\usepackage{amsfonts}

\usepackage{bm}
\usepackage{enumerate} 
\usepackage{amssymb,amsmath}
\usepackage{caption}
\usepackage{subcaption}
\usepackage{accents} 

\usepackage{dsfont}

\usepackage{stackengine}

\usepackage{accents} 

\let\Horig\H

\usepackage{tikz}
\usetikzlibrary{automata,topaths}
\usetikzlibrary{shapes}
\usetikzlibrary{plotmarks}
 
\usepackage{float} 

\usepackage{color} 
\definecolor{lightblue}{rgb}{0,0.2,0.5}
\usepackage[colorlinks=true, urlcolor=lightblue,linkcolor=lightblue, citecolor=lightblue]{hyperref}

\DeclareMathAlphabet{\eufrak}{U}{}{}{} 
\SetMathAlphabet\eufrak{normal}{U}{euf}{m}{n}
\SetMathAlphabet\eufrak{bold}{U}{euf}{b}{n}

\newcommand{\R}{\mathbb{R}}

\newcommand{\E}{\mathbb{E}}
\newcommand{\IP}{\mathbb{P}}

\newcommand{\bone}{{\bf 1}}

\newcommand{\N}{\mathbb{N}}

\newtheorem{prop}{Proposition}[section]
\newtheorem{assumption}[prop]{Assumption}
\newtheorem{lemma}[prop]{Lemma}
\newtheorem{definition}[prop]{Definition}
\newtheorem{corollary}[prop]{Corollary}

\newtheorem{example}[prop]{Example}

\makeatletter
\newcommand*\rel@kern[1]{\kern#1\dimexpr\macc@kerna}
\newcommand*\widebar[1]{
  \begingroup
  \def\mathaccent##1##2{
    \rel@kern{0.8}
    \overline{\rel@kern{-0.8}\macc@nucleus\rel@kern{0.2}}
    \rel@kern{-0.2}
  }
  \macc@depth\@ne
  \let\math@bgroup\@empty \let\math@egroup\macc@set@skewchar
  \mathsurround\z@ \frozen@everymath{\mathgroup\macc@group\relax}
  \macc@set@skewchar\relax
  \let\mathaccentV\macc@nested@a
  \macc@nested@a\relax111{#1}
  \endgroup
}
\makeatother

\makeatletter
\DeclareRobustCommand\widecheck[1]{{\mathpalette\@widecheck{#1}}}
\def\@widecheck#1#2{
    \setbox\z@\hbox{\m@th$#1#2$}
    \setbox\tw@\hbox{\m@th$#1
       \widehat{
          \vrule\@width\z@\@height\ht\z@
          \vrule\@height\z@\@width\wd\z@}$}
    \dp\tw@-\ht\z@
    \@tempdima\ht\z@ \advance\@tempdima2\ht\tw@ \divide\@tempdima\thr@@
    \setbox\tw@\hbox{
       \raise\@tempdima\hbox{\scalebox{1}[-1]{\lower\@tempdima\box
\tw@}}}
    {\ooalign{\box\tw@ \cr \box\z@}}}
\makeatother

\def\({\left(}
\def\){\right)}

\def\[{\left[}
\def\]{\right]}
\def\real{{\mathord{\mathbb R}}}
\def\N{{\mathord{\mathbb N}}}

\def\P{\mathbb{P}}

\newenvironment{Proof}{\removelastskip\par\medskip
\noindent{\em Proof.} \rm}{\penalty-20\null\hfill$\square$\par\medbreak}

\allowdisplaybreaks

\numberwithin{equation}{section}

\usepackage{graphicx}
\usepackage{flushend,cuted}
\usepackage{bm}
\usepackage{tabularx}

\usepackage{indentfirst}
\usepackage{amssymb}
\usepackage{xparse}
\usepackage{tikz}
\usepackage{mdwlist}
\usepackage{tkz-graph}
\usepackage{textgreek} 

\GraphInit[vstyle = Shade]
\usetikzlibrary[intersections,
positioning,
petri,
backgrounds,
fit,
decorations.pathmorphing,
arrows,
arrows.meta,
bending,
calc,
intersections,
through,
backgrounds,
shapes.geometric,
quotes,
matrix,
trees,
shapes.symbols,
graphs,
math,
patterns,
external,
scopes,
matrix,
lindenmayersystems,
shapes.callouts,
shapes.misc,
angles,
shapes.arrows,
shadings]

\usetikzlibrary{arrows,automata,shapes,positioning,decorations.pathmorphing,snakes}

\tikzset{snake it/.style={-stealth,
decoration={snake, 
    amplitude = .4mm,
    segment length = 2mm,
    post length=0.9mm},decorate}}

\usetikzlibrary{matrix,calc}

\newcommand{\convexpath}[2]{
[   
    create hullnodes/.code={
        \global\edef\namelist{#1}
        \foreach [count=\counter] \nodename in \namelist {
            \global\edef\numberofnodes{\counter}
            \node at (\nodename) [draw=none,name=hullnode\counter] {};
        }
        \node at (hullnode\numberofnodes) [name=hullnode0,draw=none] {};
        \pgfmathtruncatemacro\lastnumber{\numberofnodes+1}
        \node at (hullnode1) [name=hullnode\lastnumber,draw=none] {};
    },
    create hullnodes
]
($(hullnode1)!#2!-90:(hullnode0)$)
\foreach [
    evaluate=\currentnode as \previousnode using \currentnode-1,
    evaluate=\currentnode as \nextnode using \currentnode+1
    ] \currentnode in {1,...,\numberofnodes} {
-- ($(hullnode\currentnode)!#2!-90:(hullnode\previousnode)$)
  let \p1 = ($(hullnode\currentnode)!#2!-90:(hullnode\previousnode) - (hullnode\currentnode)$),
    \n1 = {atan2(\y1,\x1)},
    \p2 = ($(hullnode\currentnode)!#2!90:(hullnode\nextnode) - (hullnode\currentnode)$),
    \n2 = {atan2(\y2,\x2)},
    \n{delta} = {-Mod(\n1-\n2,360)}
  in 
    {arc [start angle=\n1, delta angle=\n{delta}, radius=#2]}
}
-- cycle
}

\tikzset{hide labels/.style={every label/.append style={text opacity=0}}}

\begin{document}
\title{
\huge
Normal approximation of subgraph counts in the random-connection model
} 

\author{
  Qingwei Liu\footnote{\href{mailto:qingwei.liu@ntu.edu.sg}{qingwei.liu@ntu.edu.sg}}
  \qquad
      Nicolas Privault\footnote{
\href{mailto:nprivault@ntu.edu.sg}{nprivault@ntu.edu.sg}
}
  \\
\small
Division of Mathematical Sciences
\\
\small
School of Physical and Mathematical Sciences
\\
\small
Nanyang Technological University
\\
\small
21 Nanyang Link, Singapore 637371
}

\maketitle

\vspace{-0.5cm}

\begin{abstract} 
 This paper derives normal approximation results for subgraph counts written as multiparameter stochastic integrals in a random-connection model based on a Poisson point process. By combinatorial arguments we express the cumulants of general subgraph counts using sums over connected partition diagrams, after cancellation of terms obtained by M\"obius inversion. Using the Statulevi\v{c}ius condition, we deduce convergence rates in the Kolmogorov distance by studying the growth of subgraph count cumulants as the intensity of the underlying Poisson point process tends to infinity. Our analysis covers general subgraphs in the dilute and full random graph regimes, and tree-like subgraphs in the sparse random graph regime.
\end{abstract}

 
\baselineskip0.7cm

\section{Introduction}
\noindent
 This paper treats the asymptotic behavior of
 random subgraph counts in the random-connection model, 
 which is used to model physical systems in e.g. wireless networks, 
 complex networks, 
 and statistical mechanics. 
 Our approach relies on the study of cumulant growth rates
 as the intensity of the underlying Poisson point process tends to infinity. 
 
\medskip 

The distributional approximation of subgraph counts has attracted
significant interest in the random graph literature. 
 In \cite{rucinski}, conditions for the asymptotic normality
 of renormalized subgraph counts have been obtained
 in the Erd{\H o}s-R\'enyi random graph model \cite{ER,G}. 
 In \cite[Chapter~3]{penrosebk}, non-quantitative central limit theorems have been obtained for subgraph and component counts
 on random geometric graphs.
   Rates of normal convergence with respect to the Kolmogorov distance
   have been obtained for certain random
   functionals on random geometric graphs in
   \cite{schulte} using Poisson $U$-statistics,
   see also \cite{lachiezerey4} for the use of stabilizing functionals.
 
\medskip 
 
 In the Erd{\H o}s-R\'enyi setting,
 those results have been made more precise in \cite{BKR} by
 the derivation of convergence rates in the Wasserstein distance via the
 Stein method. 
 They have also been strengthened 
 in the Kolmogorov distance for triangle counts
 in \cite{reichenbachsAoP}, 
 and for general subgraph counts in \cite{PS2}.
 The case of triangles has also been treated in \cite{roellin2}
 by the {S}tein-{T}ikhomirov method, which has been extended 
 to general subgraphs in \cite{rednos}.    
 In \cite{khorunzhiy}, the counts of line ($X$-model) and cycles ($Y$-model) 
 in discrete Erd{\H o}s-R\'enyi models 
 have been analyzed via the asymptotic behavior of their cumulants.
 
\medskip 
  
 The random-connection model is a natural generalization of the 
 Erd\H os-R\'enyi random graph in which vertices are randomly located
 and can be connected
 with location-dependent probabilities
 $H(x,y)\in [0,1]$.
 Obtaining normal approximation error bounds 
 in the random-connection model
 with a general $[0,1]$-valued random connection function 
 is more difficult due to the additional layer of complexity coming from
 the randomness of vertex locations.
 
\medskip 
   
 Regarding convergence rates, in \cite{LNS21},
 a central limit theorem and Berry-Esseen convergence rates have
 been presented and applied to the number of components isomorphic
 to a given finite connected graph in the random-connection model,
 together with a study of ﬁrst moments and covariances.
 In \cite{zhangzs}, Berry-Esseen convergence rates have
 been obtained for subgraph counts
 in the binomial random-connection (graphon) model. 
 However, those results do not cover the case of
 general subgraph counting in the Poisson random-connection model. 

\medskip 
 
    The Malliavin-Stein machinery on Poisson space \cite{lastpeccatipenrose} has been applied to in \cite{LNS21} the numbers of components isomorphic to a given graph in the random-connection model, using edge marked Poisson processes.
   Recently, a central limit theorem has been derived in \cite{can2022} for the counts of induced subgraphs in the random-connection model
   using the edge-marking structure of \cite{LNS21}, 
   under a weak stabilizing condition originating from \cite{penrose01}. 
   However, as pointed out in Remark~2.5-$(i)$ of \cite{can2022},
   no convergence rates are derived by this method, as 
   the strong stabilization condition of \cite{penrose05,lachiezerey4}
   is not satisfied by general functionals
   when the connection function $H(x,y)$ is $(0,1)$-valued. 
   On the other hand, the cumulant method,
   combined with the use of partition diagrams,
   enables us to establish a quantitative central limit theorem for functionals
   in the random-connection model. 

\medskip 

 In this paper, we derive normal approximation rates under a mild condition
 on the connection function $H(x,y)$ of the random-connection model,
 by deriving growth rates of cumulants written as sums over
 connected partitions, see Propositions~\ref{t1} and \ref{th6.4}.
 {Related cumulant bounds have been obtained
   in the Erd{\H o}s-R\'enyi model,
   cf. Proposition~10.1.2 in \cite{feray}. 
   However, }
 to the best of our knowledge, this is the first time that the normal approximation
 of subgraph counts with convergence rates is established
 in the random-connection model.
 
\medskip 

 In comparison with \cite{khorunzhiy}, which also uses the cumulant method, 
 we obtain convergence rates in the Kolmogorov distance and our results are 
 not restricted to line and cycle graphs, as they cover more general subgraphs,
 see Corollaries~\ref{c01}-\ref{c01-2}. 
 Furthermore, various random graph regimes are discussed.
 In addition,
   we show in Section~\ref{rgg} that our approach can be
   specialized to derive Kolmogorov rates for
   subgraph counting in the setting of random geometric graphs,
   see Corollary~\ref{jdkj10}. 
 
\medskip 

 A number of probabilistic conclusions can be derived
 from the behavior of cumulants of random variables using the 
Statulevi\v{c}ius condition, 
 including convergence rates in the Kolmogorov distance
 and moderate deviation principles, see
 \cite{saulis},
 \cite{doring},
 \cite{doering}. 
{
 In stochastic geometry, 
 the cumulant method has also been applied
  to Poisson cylinder processes \cite{heinrich}, 
  and to the volumes of simplices
  in Poisson-Delaunay tessellations \cite{gusakova}, 
  to the Boolean model \cite{heinrich2}, 
  and to random $m$-dependent ﬁelds \cite{hipp}.
}
 In \cite{grotethale18,thale18}, this method has been used to 
 derive concentration inequalities,
 normal approximation with error bounds,
 and moderate deviation principles for random polytopes. 
   
\medskip 
 
Given $\mu$ a 
 {finite} diffuse measure on $\R^d$,
 we consider a random-connection model
 based on an underlying Poisson point process $\Xi$ on $\R^d$
 with intensity of the form $\lambda\mu(\mathrm{d}x)$, in which 
 any two vertices $x,y$ in $\Xi$ are connected
 with the probability $H_\lambda(x,y):= c_\lambda H(x,y) \in [0,1]$,
 where $H_\lambda$ is the connection function of the model. 
 Here, we investigate the limiting behavior
 of the count $N_G$ of a given subgraph $G$
 as the intensity $\lambda$ of the underlying Poisson point process on $\R^d$
 tends to infinity. 
 To this end, we use the combinatorics of the cumulants $\kappa_n(N_G)$ 
 based on moment expressions obtained in \cite{prkhp} for
 multiparameter stochastic integrals in the random-connection model.
 
 \medskip 

Using partition diagrams and dependency graph arguments,
we start by showing in Proposition~\ref{mainthm-1}
that the (virtual) cumulants of a random functional admitting
a certain connectedness factorization property \eqref{dia-factoriz}  
can be expressed as sums over connected partition diagrams, 
 generalizing Lemma~2 in \cite{MalyshevMinlos91}.  
 A related result has been obtained in \cite{jansen}
 in the particular case of two-parameter Poisson stochastic
 integrals, in relation to 
 cluster expansions for Gibbs point processes in statistical
 mechanics. 
 In Proposition~\ref{p01-1}, we apply Proposition~\ref{mainthm-1} 
 to express the cumulants of multiparameter stochastic integrals,
 for which this factorization property can be checked from 
 the moment formulas for 
 multiparameter stochastic integrals computed in Proposition~\ref{p01-1-0}. 
 
\medskip 
 
 Such expressions allow us to determine the dominant terms in the growth of
 cumulants as the intensity $\lambda$ of the underlying point process tends to infinity,
 by estimating the counts of vertices and edges in connected partition diagrams
 as in \cite{khorunzhiy}. 
 We work under a mild condition \eqref{integ-connecting3} 
 which is satisfied by e.g. any translation-invariant
 continuous connection function $H : \real^d\times \real^d \to [0,1]$ non vanishing at $0$, such as the {Rayleigh} connection function
  given by $H(x,y) = e^{ - \beta \Vert x - y\Vert^2}$, $x,y\in \real^d$, for some $\beta > 0$. 
 
\medskip 
 
 For our analysis of cumulant behavior
 we identify the leading terms in the sum \eqref{cumulant-diagram1}
 over connected partition diagrams. 
 When $G$ is a connected graph with $|V(G)|=r$ vertices,
 satisfying
 Assumption~\ref{a61} in the dilute regime \eqref{fjnldsf}
 with $\lambda^{-1/\zeta } \ll c_\lambda \leq K$,
 where $\zeta \geq 1$ is defined in \eqref{fjkldf}, 
 the dominant terms are given by connected partition diagrams with the
 highest number of blocks, 
 see also \cite{privaultkhops}
 in the case of $k$-hop counting on the line.
  In Proposition~\ref{t1} this yields the cumulant bounds 
$$
 (n-1)! c_\lambda^{n |E(G)| } ( K_1 \lambda )^{1+(r-1)n} 
 \leq 
  \kappa_n(N_G)
\leq 
n!^r c_\lambda^{n |E(G) |} ( K_2 \lambda )^{1+(r-1)n},
\quad \lambda \geq 1, 
$$
for some constants $K_1$, $K_2>0$ independent of $\lambda, n\geq 1$,
where $E(G)$ denotes the set of edges of $G$.
 From the {Statulevi\v{c}ius condition}
 \eqref{Statuleviciuscond2} below, see \cite{rudzkis,doering},
 letting $\Phi$ denote the cumulative distribution function of the standard normal distribution, 
 we deduce the Kolmogorov distance bound 
$$
\sup_{x\in \real}
\big| \P \big( \widetilde{N}_G \leq x \big) - \Phi (x) \big| \leq
\frac{C }{\lambda^{1/(4r - 2)}},
\qquad \lambda \to \infty, 
$$ 
 for the normalized subgraph count $\widetilde{N}_G$, 
see Corollary~\ref{c01}, and a moderate deviation principle, 
 see Corollary~\ref{c01-2-0}. 
  
\medskip

 In the sparse regime \eqref{fjnldsf-2} where
 $ c_\lambda \leq \lambda^{-\alpha}$ for some $\alpha \geq 1$,
 the maximal rate $\lambda^{
 \alpha         -(\alpha - 1)r 
          }$
 is attained for $G$ a tree-like graph, and
 in Proposition~\ref{th6.4} we obtain the cumulant bounds 
$$ 
    \nonumber 
            (K_1)^r
    \lambda^{
 \alpha     -(\alpha - 1)r 
      }
     \leq 
  \kappa_n(N_G)
    \leq
  n!^r
    (K_2)^r
  \lambda^{
   \alpha -(\alpha - 1)r 
    } 
  , \quad \lambda \geq 1, 
$$ 
  if $G$ is a tree, and
$$ 
  \nonumber 
    (K_1)^r 
  \lambda^{r-\alpha |E(G)|}
  \leq 
  \kappa_n(N_G)
  \leq
    n!^r
      (K_2)^r
    \lambda^{r-\alpha |E(G)|}, \quad \lambda \geq 1, 
$$ 
  if $G$ is a not a tree, such as e.g. a cycle graph.
  As a consequence of the {Statulevi\v{c}ius condition}
  \eqref{Statuleviciuscond2}, 
 when $G$ is a tree
 we find the Kolmogorov distance bound 
$$  
\sup_{x\in \real}
\big| \P \big( \widetilde{N}_G \leq x \big) - \Phi (x) \big| \leq
 C \lambda^{
    - (
\alpha    -(\alpha - 1)r 
        ) / ( 4r - 2) }
, \qquad \lambda \to \infty, 
$$ 
 provided that $1 \leq \alpha < r/(r-1)$, 
 see Corollary~\ref{c01-2}. 
 
 \medskip 
 
 Convergence rates in the Kolmogorov distances may be improved
 into classical Berry-Esseen rates when 
 the connection function $H(x,y)$ is $\{0,1\}$-valued, e.g. in disk models
 as in \cite{privaultkhops}, 
 by representing subgraph counts as multiple Poisson stochastic integrals
 and using the fourth moment theorem for $U$-statistics and sums of
 multiple stochastic integrals Corollary~4.10 in \cite{eichelsbacher}, 
 see also Theorem~3 in \cite{lachieze-rey} 
 or Theorem~6.3 in \cite{PS4} for Hoeffding decompositions.
 On the other hand, the study of stabilizing functionals \cite{penrose05,lachiezerey4} yields normal approximation with rates for random functionals
 represented as sums of stabilizing score functions on random geometric graphs.
 In the general case where $H(x,y)$ is $[0,1]$-valued,
 both methods no longer apply, which is why we rely on the
 {cumulant method}
 and the
 {Statulevi\v{c}ius condition},
 which in turn may yield suboptimal convergence rates. 
{
 Recently, moderate deviation principles
 have been obtained by the cumulant method
 for functionals of Poisson point processes 
 in \cite{schulte-thaele}, with application
 to subgraph counting in random geometric graphs.} 
{However, \cite{schulte-thaele}
   does not cover
  the random-connection model with a general connection
  function
  $H(x,y) \in [0,1]$.}

 \medskip

 This paper is organized as follows.
 Sections~\ref{s2} and \ref{s3} introduce the preliminary
 framework and notations on connected partition diagrams
 and combinatorics of virtual cumulants that will be used for 
 the expression of cumulants of multiparameter stochastic integrals
 in Section~\ref{s4} and for subgraph counts in Section~\ref{s5}. 
 Those expressions are applied in Section~\ref{s6}
 to derive cumulant growth rates in the random-connection model, with application to Kolmogorov rates in subgraph counting via the {Statulevi\v{c}ius condition} in Section~\ref{s6-1}.
 In Section~\ref{rgg}, normal approximation for subgraph counts on the random geometric graph is discussed under different limiting regimes. 

\subsubsection*{Preliminaries} 
 \noindent
 Consider a Poisson point process $\Xi$ on $\R^d$, $d \geq 1$, with
{$\sigma$-finite}
intensity measure $\Lambda$ on $\real^d$, 
 constructed on the space $$
 \Omega = \big\{
 \omega = \{ x_i \}_{i\in I} \subset \R^d \ : \
 \#( A \cap \omega ) < \infty 
 \mbox{ for all compact } A\in {\cal B} (\R^d) 
 \big\}
 $$
 of locally finite configurations on $\R^d$, whose elements 
 $\omega \in \Omega$ are identified with the Radon point measures 
 $\displaystyle \omega = \sum_{x\in \omega} \epsilon_x$, 
 where $\epsilon_x$ denotes the Dirac measure at $x\in \R^d$. 
 {As in} 
 \cite[Corollary~6.5]{LastPenrose17}, almost every
 element $ \omega$ of $\Omega$ 
 can be represented as $\omega =\{V_i\}_{1\leq i\leq N}$,
 where $(V_i)_{i\geq 1}$ is a random sequence 
 in $\R^d$ and a $\N\cup\{\infty\}$-valued random variable $N$.
 
\medskip

In what follows, we let $[n]:=\{1,2,\dots,n\}$ for $n\geq 1$.
In the next proposition, see Proposition~2 in \cite{prkhp}, 
which relies on Proposition~3.1 of \cite{momentpoi} 
and Lemma~2.1 of \cite{bogdan}, we express the moments of \eqref{e1}
using sums over the set
 $\Pi ( \eta \times [r])$ 
 of all partitions of the set 
 $$
 \eta \times [r] := 
 \big\{ (k,l) \ : \
 k\in \eta, \ l = 1,\ldots , r \big\},
 \quad
 n,r\geq 1, \ \eta \subset [n].
 $$
  
 \begin{prop}
\label{p01-1-0} 
 Given $r\geq 2$, consider a 
 connected graph $G$ with $r$ vertices, 
 and a {bounded} 
 measurable process of the form 
$$
 u (x_1,\ldots , x_r ) := \prod_{\{i,j\} \in E(G)} v(x_i,x_j), 
$$
 where $v(x,y)$ is a {bounded}
 random process $v(x,y)$
 independent of the underlying Poisson point process~$\Xi$. 
 Then, the $n$-$th$ moment of the multiparameter stochastic integral
\begin{equation}
\label{e1}
 \sum_{\{ V_1,\dots,V_r\} \subset \omega } \ u (V_1,\ldots , V_r ) 
 = \int_{(\real^d)^r} u (x_1,\ldots , x_r ) \omega (\mathrm{d}x_1)\cdots
 \omega (\mathrm{d}x_r),
 \quad n \geq 1, 
\end{equation} 
is given by the summation 
\begin{equation} 
\label{nthexpectation0-0} 
\sum_{\rho \in \Pi ( [n] \times [r] )}
\int_{(\R^d )^{|\rho |}}
\E \left[
  \prod_{k=1}^n
\prod_{
\{i,j\}\in E(G)}
v \big(x^\rho_{k,i},x^\rho_{k,j}\big) 
\right] 
   \ \! \prod_{\eta \in \rho }
   \Lambda(\mathrm{d}x_\eta ),
\end{equation}
 where we let $x_{k,l}^\rho:=x_\eta$ whenever $(k,l)\in \eta$,
 for $\rho \in \Pi([n]\times [r])$ and $\eta\in\rho$.
\end{prop}
\section{Set partitions and diagram connectivity} 
\label{s2}
Given $\eta$ a finite set, we denote by $\Pi ( \eta )$ the collection
of its set partitions, and we let $|\sigma|$ denote the number of blocks in any partition $\sigma \in \Pi ( \eta )$. 
Given $\rho,\sigma$ two set partitions, we say that $\sigma$ is coarser than $\rho$,
or that $\rho$ is finer than $\sigma$,
and we write $\rho\preceq\sigma$,
if every block in $\sigma$ is a combination of blocks in $\rho$. 
We also denote by $\rho\vee\sigma$ the finest partition which is coarser than $\rho$ and $\sigma$, and by $\rho\wedge\sigma$ the coarsest partition that is finer than $\rho$ and $\sigma$.
We let $\widehat{0}$ be the finest partition, which is made of a single element in each block, and we let $\widehat{1}$ be the coarsest (one-block) partition.
In general, given any graph $G$ we denote by $V(G)$ the set of its vertices, and by $E(G)$ the set of its edges.

\medskip

Our study of cumulants and moments of functionals of random fields
relies on partition diagrams, see \cite{MalyshevMinlos91,khorunzhiy,peccatitaqqu} 
and references therein for additional background. 
 In the sequel, for $n,r\geq 1$ 
 we let 
$\pi_\eta : = (\pi_i)_{i\in \eta} \in \Pi ( \eta \times [r])$ denote the
 partition made of the $|\eta|$ blocks
 of size $r$ given by 
 $\pi_k := \{ (k,1), \ldots , (k,r) \}$, $k\in \eta$. 
\begin{definition}
  Let $n,r \geq 1$.
  Given $\eta\subset [n]$ and $\rho\in\Pi(\eta\times[r])$
 a partition of $\eta \times [r]$, we denote by $\Gamma (\rho,\pi_\eta )$
 the diagram, or graphical representation of the partition $\rho$,
 constructed by: 
\begin{enumerate}[\rm 1.] 
\item arranging the elements of $\eta \times [r]$
 into an 
 array of $|\eta |$ rows and $r$ columns, and
\item
 connecting all elements within a same block of $\rho$ 
 by a tree graph. 
\end{enumerate} 
In addition, we say that the partition diagram $\Gamma(\rho,\pi_\eta )$
 is connected when $\rho\vee\pi_\eta=\widehat{1}$. 
\end{definition}
\noindent 
 For shortness of notation, in the sequel  
 we say that a partition $\rho$ is connected 
 when its diagram $\Gamma ( \rho , \pi )$ is connected. 
 For example, taking $\eta := \{2,3,5,8,10\}$, given the partitions 
\begin{align*}
  \rho = \big\{
 & \{(2,1),(3,1),(3,2),(3,3)\}, 
\{(2,2),(2,3),(2,4),(3,4)\},
\{(5,1)\},
\{(5,2),(8,2)\},
\\
& \{(5,3)\},
\{(5,4),(8,3)\},
\{(8,1),(10,1)\},
\{(8,4)\},
\{(10,2),(10,3),(10,4)\}\big\}
\end{align*} 
and
\begin{align*}
  \sigma = \big\{ & 
  \{(2,1),(3,1)\},
  \{(2,2)\},
  \{(2,3),(3,4)\},
  \{(2,4)\},
  \{(3,2),(5,2),(8,2)\},
  \\
  &
  \{(3,3),(5,4),(8,3),(10,2)\},
  \{(5,1)\},
  \{(5,3)\},
  \{(8,1),(10,1)\},
  \{(8,4)\},
  \{(10,3)\},
  \{(10,4)\}
  \big\}, 
\end{align*} 
of $\eta \times [4]$,
Figure~\ref{fig:diagram0}-$a)$ presents an example of a non-connected partition diagram  
$\Gamma ( \rho , \pi)$,
and Figure~\ref{fig:diagram0}-$b)$ presents an example of a connected partition diagram $\Gamma ( \sigma , \pi)$. 

\begin{figure}[H]
\captionsetup[subfigure]{font=footnotesize}
\centering
\subcaptionbox{Non-connected partition diagram $\Gamma(\rho,\pi)$.}[.49\textwidth]{
\begin{tikzpicture}[scale=0.9] 
\draw[black, thick] (0,0) rectangle (5,6);
\draw[very thick,dashed,red] (0.3,3.5) -- (4.7,3.5);
\node[anchor=east,font=\small] at (0.8,5) {2};
\node[anchor=east,font=\small] at (0.8,4) {3};
\node[anchor=east,font=\small] at (0.8,3) {5};
\node[anchor=east,font=\small] at (0.8,2) {8};
\node[anchor=east,font=\small] at (0.8,1) {10};

\node[anchor=south,font=\small] at (1,0) {1};
\node[anchor=south,font=\small] at (2,0) {2};
\node[anchor=south,font=\small] at (3,0) {3};
\node[anchor=south,font=\small] at (4,0) {4};

\filldraw [gray] (1,1) circle (2pt);
\filldraw [gray] (2,1) circle (2pt);
\filldraw [gray] (3,1) circle (2pt);
\filldraw [gray] (4,1) circle (2pt);
\filldraw [gray] (1,2) circle (2pt);
\filldraw [gray] (2,2) circle (2pt);
\filldraw [gray] (3,2) circle (2pt);
\filldraw [gray] (4,2) circle (2pt);
\filldraw [gray] (1,3) circle (2pt);
\filldraw [gray] (2,3) circle (2pt);
\filldraw [gray] (3,3) circle (2pt);
\filldraw [gray] (4,3) circle (2pt);
\filldraw [gray] (2,3) circle (2pt);
\filldraw [gray] (1,4) circle (2pt);
\filldraw [gray] (2,4) circle (2pt);
\filldraw [gray] (3,4) circle (2pt);
\filldraw [gray] (4,4) circle (2pt);
\filldraw [gray] (1,5) circle (2pt);
\filldraw [gray] (2,5) circle (2pt);
\filldraw [gray] (3,5) circle (2pt);
\filldraw [gray] (4,5) circle (2pt);

\draw[very thick] (1,5) -- (1,4) -- (2,4) -- (3,4);
\draw[very thick] (2,5) -- (3,5) -- (4,5) -- (4,4);

\draw[very thick] (1,2) -- (1,1);
\draw[very thick] (2,3) -- (2,2);
\draw[very thick] (2,1) -- (3,1) -- (4,1);
\draw[very thick] (3,2) -- (4,3);

  \begin{pgfonlayer}{background}
    \filldraw [line width=4mm,black!3]
      (0.2,0.2)  rectangle (4.8,5.8);
  \end{pgfonlayer}
\end{tikzpicture}}
\subcaptionbox{Connected partition diagram $\Gamma(\sigma,\pi)$.}[.49\textwidth]{
\begin{tikzpicture}[scale=0.9] 
\draw[black, thick] (0,0) rectangle (5,6);

\node[anchor=east,font=\small] at (0.8,5) {2};
\node[anchor=east,font=\small] at (0.8,4) {3};
\node[anchor=east,font=\small] at (0.8,3) {5};
\node[anchor=east,font=\small] at (0.8,2) {8};
\node[anchor=east,font=\small] at (0.8,1) {10};

\node[anchor=south,font=\small] at (1,0) {1};
\node[anchor=south,font=\small] at (2,0) {2};
\node[anchor=south,font=\small] at (3,0) {3};
\node[anchor=south,font=\small] at (4,0) {4};

\filldraw [gray] (1,1) circle (2pt);
\filldraw [gray] (2,1) circle (2pt);
\filldraw [gray] (3,1) circle (2pt);
\filldraw [gray] (4,1) circle (2pt);
\filldraw [gray] (1,2) circle (2pt);
\filldraw [gray] (2,2) circle (2pt);
\filldraw [gray] (3,2) circle (2pt);
\filldraw [gray] (4,2) circle (2pt);
\filldraw [gray] (1,3) circle (2pt);
\filldraw [gray] (2,3) circle (2pt);
\filldraw [gray] (3,3) circle (2pt);
\filldraw [gray] (4,3) circle (2pt);
\filldraw [gray] (2,3) circle (2pt);
\filldraw [gray] (1,4) circle (2pt);
\filldraw [gray] (2,4) circle (2pt);
\filldraw [gray] (3,4) circle (2pt);
\filldraw [gray] (4,4) circle (2pt);
\filldraw [gray] (1,5) circle (2pt);
\filldraw [gray] (2,5) circle (2pt);
\filldraw [gray] (3,5) circle (2pt);
\filldraw [gray] (4,5) circle (2pt);

\draw[very thick] (1,5) -- (1,4); 
\draw[very thick] (3,5) -- (4,4);

\draw[very thick] (1,2) -- (1,1);
\draw[very thick] (2,2) -- (2,4);
\draw[very thick] (2,1) -- (3,2) -- (4,3) -- (3,4);

  \begin{pgfonlayer}{background}
    \filldraw [line width=4mm,black!3]
      (0.2,0.2)  rectangle (4.8,5.8);
  \end{pgfonlayer}
\end{tikzpicture}}
\caption{Two examples of partition diagrams with $\eta = \{2,3,5,8,10\}$, $n=10$, $r=4$.}
\label{fig:diagram0}
\end{figure}

\vspace{-.4cm}

\noindent 
Note that the above notion of connected partition diagram is distinct from
that of irreducible partition, see, e.g., \cite{eabender}.
\begin{definition}
  Let $n\geq 1$,
  $G$ a connected graph with $|V(G)| = r$ vertices, $r \geq 1$,  
  and consider $G_1,\ldots , G_n$ copies of $G$ respectively built on
  $\pi_1,\ldots , \pi_n$.
  Let also $\rho \in\Pi( [n] \times[r])$
  be a partition of $ [n] \times[r]$. 
\begin{enumerate}[\rm 1.]
\item
  We let $\widetilde{\rho}_G$ be the multigraph constructed
  on the blocks of $\rho$ 
  by adding an edge between two blocks $\rho_1,\rho_2$ of the
  partition $\rho$ whenever there exist $(k,l_1)\in \rho_1$
  and $(k,l_2)\in \rho_2$ such that $(l_1,l_2)$ is an edge in $G_k$.
  \item 
  We let $\rho_G$ be the graph constructed
  on the blocks of $\rho$
  by removing redundant edges
  in $\widetilde{\rho}_G$,
  so that at most one edge remains between any two blocks $\rho_1,\rho_2\in\rho$. 
 \end{enumerate}
\end{definition}

\begin{figure}[H]
\captionsetup[subfigure]{font=footnotesize}
\centering
\subcaptionbox{Diagram $\Gamma(\rho,\pi)$ and multigraph $\widetilde{\rho}_G$ in blue.}[.49\textwidth]{
\begin{tikzpicture}[scale=0.9] 
\draw[black, thick] (0,0) rectangle (5,6);

\node[anchor=east,font=\small] at (0.8,5) {1};
\node[anchor=east,font=\small] at (0.8,4) {2};
\node[anchor=east,font=\small] at (0.8,3) {3};
\node[anchor=east,font=\small] at (0.8,2) {4};
\node[anchor=east,font=\small] at (0.8,1) {5};

\node[anchor=south,font=\small] at (1,0) {1};
\node[anchor=south,font=\small] at (2,0) {2};
\node[anchor=south,font=\small] at (3,0) {3};
\node[anchor=south,font=\small] at (4,0) {4};

\filldraw [gray] (1,1) circle (2pt);
\filldraw [gray] (2,1) circle (2pt);
\filldraw [gray] (3,1) circle (2pt);
\filldraw [gray] (4,1) circle (2pt);
\filldraw [gray] (1,2) circle (2pt);
\filldraw [gray] (2,2) circle (2pt);
\filldraw [gray] (3,2) circle (2pt);
\filldraw [gray] (4,2) circle (2pt);
\filldraw [gray] (1,3) circle (2pt);
\filldraw [gray] (2,3) circle (2pt);
\filldraw [gray] (3,3) circle (2pt);
\filldraw [gray] (4,3) circle (2pt);
\filldraw [gray] (2,3) circle (2pt);
\filldraw [gray] (1,4) circle (2pt);
\filldraw [gray] (2,4) circle (2pt);
\filldraw [gray] (3,4) circle (2pt);
\filldraw [gray] (4,4) circle (2pt);
\filldraw [gray] (1,5) circle (2pt);
\filldraw [gray] (2,5) circle (2pt);
\filldraw [gray] (3,5) circle (2pt);
\filldraw [gray] (4,5) circle (2pt);

\draw[very thick] (1,5) -- (1,4) -- (2,4) -- (3,4);
\draw[very thick] (2,5) -- (3,5) -- (4,5) -- (4,4);

\draw[very thick] (1,2) -- (1,1);
\draw[very thick] (2,3) -- (2,2);
\draw[very thick] (2,1) -- (3,1) -- (4,1);
\draw[very thick] (3,2) -- (4,3);

\draw[thick,dash dot,blue] (1,5) .. controls (1.5,5+.5) .. (2,5);

\draw[thick,dash dot,blue] (3,4) .. controls (3.5,4+.5) .. (4,4);
\draw[thick,dash dot,blue] (2,4) .. controls (3,4-.5) .. (4,4);
\foreach \i in {2,...,3}
         {
\draw[thick,dash dot,blue] (1,\i) .. controls (1.5,\i+.5) .. (2,\i);
\draw[thick,dash dot,blue] (3,\i) .. controls (3.5,\i+.5) .. (4,\i);
\draw[thick,dash dot,blue] (2,\i) .. controls (3,\i-.5) .. (4,\i);
} 

\draw[thick,dash dot,blue] (1,1) .. controls (1.5,1+.5) .. (2,1);

  \begin{pgfonlayer}{background}
    \filldraw [line width=4mm,black!3]
      (0.2,0.2)  rectangle (4.8,5.8);
  \end{pgfonlayer}
\end{tikzpicture}}
\subcaptionbox{Diagram $\Gamma(\rho ,\pi)$ and graph $\rho_G$ in red.}[.49\textwidth]{
\begin{tikzpicture}[scale=0.9] 
\draw[black, thick] (0,0) rectangle (5,6);

\node[anchor=east,font=\small] at (0.8,5) {1};
\node[anchor=east,font=\small] at (0.8,4) {2};
\node[anchor=east,font=\small] at (0.8,3) {3};
\node[anchor=east,font=\small] at (0.8,2) {4};
\node[anchor=east,font=\small] at (0.8,1) {5};

\node[anchor=south,font=\small] at (1,0) {1};
\node[anchor=south,font=\small] at (2,0) {2};
\node[anchor=south,font=\small] at (3,0) {3};
\node[anchor=south,font=\small] at (4,0) {4};

\filldraw [gray] (1,1) circle (2pt);
\filldraw [gray] (2,1) circle (2pt);
\filldraw [gray] (3,1) circle (2pt);
\filldraw [gray] (4,1) circle (2pt);
\filldraw [gray] (1,2) circle (2pt);
\filldraw [gray] (2,2) circle (2pt);
\filldraw [gray] (3,2) circle (2pt);
\filldraw [gray] (4,2) circle (2pt);
\filldraw [gray] (1,3) circle (2pt);
\filldraw [gray] (2,3) circle (2pt);
\filldraw [gray] (3,3) circle (2pt);
\filldraw [gray] (4,3) circle (2pt);
\filldraw [gray] (2,3) circle (2pt);
\filldraw [gray] (1,4) circle (2pt);
\filldraw [gray] (2,4) circle (2pt);
\filldraw [gray] (3,4) circle (2pt);
\filldraw [gray] (4,4) circle (2pt);
\filldraw [gray] (1,5) circle (2pt);
\filldraw [gray] (2,5) circle (2pt);
\filldraw [gray] (3,5) circle (2pt);
\filldraw [gray] (4,5) circle (2pt);

\draw[very thick] (1,5) -- (1,4) -- (2,4) -- (3,4);
\draw[very thick] (2,5) -- (3,5) -- (4,5) -- (4,4);

\draw[very thick] (1,2) -- (1,1);
\draw[very thick] (2,3) -- (2,2);
\draw[very thick] (2,1) -- (3,1) -- (4,1);
\draw[very thick] (3,2) -- (4,3);

\draw[thick,dash dot,purple] (1,5) .. controls (1.5,5.5) .. (2,5);

\draw[thick,dash dot,purple] (1,3) .. controls (1.5,3.5) .. (2,3);
\draw[thick,dash dot,purple] (2,3) .. controls (3,2.5) .. (4,3);
\draw[thick,dash dot,purple] (3,3) .. controls (3.5,3.5) .. (4,3);

\draw[thick,dash dot,purple] (1,2) .. controls (1.5,2.5) .. (2,2);
\draw[thick,dash dot,purple] (2,2) .. controls (3,1.5) .. (4,2);
\draw[thick,dash dot,purple] (3,2) .. controls (3.5,2.5) .. (4,2);

\draw[thick,dash dot,purple] (1,1) .. controls (1.5,1.5) .. (2,1);

  \begin{pgfonlayer}{background}
    \filldraw [line width=4mm,black!3]
      (0.2,0.2)  rectangle (4.8,5.8);
  \end{pgfonlayer}
\end{tikzpicture}}
\caption{Diagram and graphs $G$, $\rho_G$, $\widetilde{\rho}_G$ with $n=5$, $r=4$.}
\label{fig:diagram1}
\end{figure}

\noindent
 Figure~\ref{fig:diagram1}-b) presents an illustration of the
multigraph $\widetilde{\rho}_G$ and graph $\rho_G$ on the blocks of $\rho$ when
$G$ is the line graph
$\{(1,2),(2,4),(3,4)\}$ on $\{1,2,3,4\}$. 

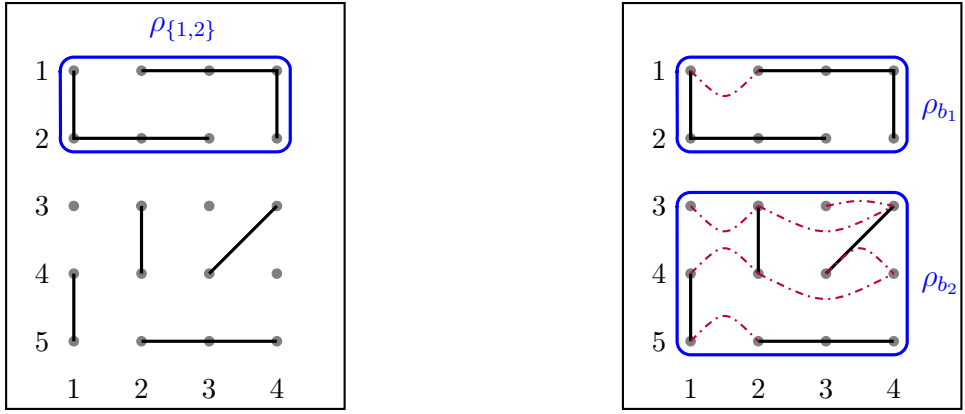
\begin{figure}[H]
\captionsetup[subfigure]{font=footnotesize}
\centering
\subcaptionbox{Diagram $\Gamma(\rho,\pi)$ and block $\rho_{\{1,2\}}$.}[.49\textwidth]{
\begin{tikzpicture}[scale=0.9,hide labels]
\tikzstyle{VertexStyle}=[shape = circle, fill = blue!20, minimum size = 0pt, scale=0., text = white, hide labels]

\draw[black, thick] (0,0) rectangle (5,6);

\node[anchor=east,font=\small] at (0.8,5) {1};
\node[anchor=east,font=\small] at (0.8,4) {2};
\node[anchor=east,font=\small] at (0.8,3) {3};
\node[anchor=east,font=\small] at (0.8,2) {4};
\node[anchor=east,font=\small] at (0.8,1) {5};

\node[anchor=south,font=\small] at (1,0) {1};
\node[anchor=south,font=\small] at (2,0) {2};
\node[anchor=south,font=\small] at (3,0) {3};
\node[anchor=south,font=\small] at (4,0) {4};

\filldraw [gray] (1,1) circle (2pt);
\filldraw [gray] (2,1) circle (2pt);
\filldraw [gray] (3,1) circle (2pt);
\filldraw [gray] (4,1) circle (2pt);
\filldraw [gray] (1,2) circle (2pt);
\filldraw [gray] (2,2) circle (2pt);
\filldraw [gray] (3,2) circle (2pt);
\filldraw [gray] (4,2) circle (2pt);
\filldraw [gray] (1,3) circle (2pt);
\filldraw [gray] (2,3) circle (2pt);
\filldraw [gray] (3,3) circle (2pt);
\filldraw [gray] (4,3) circle (2pt);
\filldraw [gray] (2,3) circle (2pt);
\filldraw [gray] (1,4) circle (2pt);
\filldraw [gray] (2,4) circle (2pt);
\filldraw [gray] (3,4) circle (2pt);
\filldraw [gray] (4,4) circle (2pt);
\filldraw [gray] (1,5) circle (2pt);
\filldraw [gray] (2,5) circle (2pt);
\filldraw [gray] (3,5) circle (2pt);
\filldraw [gray] (4,5) circle (2pt);

\draw[very thick] (1,5) -- (1,4) -- (2,4) -- (3,4);
\draw[very thick] (2,5) -- (3,5) -- (4,5) -- (4,4);

\draw[very thick] (1,2) -- (1,1);
\draw[very thick] (2,3) -- (2,2);
\draw[very thick] (2,1) -- (3,1) -- (4,1);
\draw[very thick] (3,2) -- (4,3);

\node (1) [label=above:{}] at (1,5) {};
\node (2) [label=above:{}] at (2,5) {};
\node (3) [label=above:{}] at (3,5) {};
\node (4) [label=above:{}] at (4,5) {};

\node (5) [label=above:{}] at (1,4) {};
\node (6) [label=above:{}] at (2,4) {};
\node (7) [label=above:{}] at (3,4) {};
\node (8) [label=above:{}] at (4,4) {};

\draw[very thick,blue] \convexpath{1,4,8,5}{.2cm};

\draw[blue,line width=1mm, ->] node[font=\fontsize{12}{0}\selectfont, right=of 2, right=0.5cm, below=-0.9cm] {$~~~~~~\rho_{\{1,2\}}$};

  \begin{pgfonlayer}{background}
    \filldraw [line width=4mm,black!3]
      (0.2,0.2)  rectangle (4.8,5.8);
  \end{pgfonlayer}
\end{tikzpicture}}
\subcaptionbox{Splitting $\{\rho_{b_1},\rho_{b_2}\}$
  of $\rho$ according to $\rho_G$.}[.49\textwidth]{
\begin{tikzpicture}[scale=0.9] 
\draw[black, thick] (0,0) rectangle (5,6);

\node[anchor=east,font=\small] at (0.8,5) {1};
\node[anchor=east,font=\small] at (0.8,4) {2};
\node[anchor=east,font=\small] at (0.8,3) {3};
\node[anchor=east,font=\small] at (0.8,2) {4};
\node[anchor=east,font=\small] at (0.8,1) {5};

\node[anchor=south,font=\small] at (1,0) {1};
\node[anchor=south,font=\small] at (2,0) {2};
\node[anchor=south,font=\small] at (3,0) {3};
\node[anchor=south,font=\small] at (4,0) {4};

\filldraw [gray] (1,1) circle (2pt);
\filldraw [gray] (2,1) circle (2pt);
\filldraw [gray] (3,1) circle (2pt);
\filldraw [gray] (4,1) circle (2pt);
\filldraw [gray] (1,2) circle (2pt);
\filldraw [gray] (2,2) circle (2pt);
\filldraw [gray] (3,2) circle (2pt);
\filldraw [gray] (4,2) circle (2pt);
\filldraw [gray] (1,3) circle (2pt);
\filldraw [gray] (2,3) circle (2pt);
\filldraw [gray] (3,3) circle (2pt);
\filldraw [gray] (4,3) circle (2pt);
\filldraw [gray] (2,3) circle (2pt);
\filldraw [gray] (1,4) circle (2pt);
\filldraw [gray] (2,4) circle (2pt);
\filldraw [gray] (3,4) circle (2pt);
\filldraw [gray] (4,4) circle (2pt);
\filldraw [gray] (1,5) circle (2pt);
\filldraw [gray] (2,5) circle (2pt);
\filldraw [gray] (3,5) circle (2pt);
\filldraw [gray] (4,5) circle (2pt);

\draw[very thick] (1,5) -- (1,4) -- (2,4) -- (3,4);
\draw[very thick] (2,5) -- (3,5) -- (4,5) -- (4,4);

\draw[very thick] (1,2) -- (1,1);
\draw[very thick] (2,3) -- (2,2);
\draw[very thick] (2,1) -- (3,1) -- (4,1);
\draw[very thick] (3,2) -- (4,3);

\draw[thick,dash dot,purple] (1,5) .. controls (1.5,4.5) .. (2,5);

\draw[thick,dash dot,purple] (1,3) .. controls (1.5,2.5) .. (2,3);
\draw[thick,dash dot,purple] (2,3) .. controls (3,2.5) .. (4,3);
\draw[thick,dash dot,purple] (3,3) .. controls (3.5,3.1) .. (4,3);

\draw[thick,dash dot,purple] (1,2) .. controls (1.5,2.5) .. (2,2);
\draw[thick,dash dot,purple] (2,2) .. controls (3,1.5) .. (4,2);
\draw[thick,dash dot,purple] (3,2) .. controls (3.5,2.5) .. (4,2);

\draw[thick,dash dot,purple] (1,1) .. controls (1.5,1.5) .. (2,1);

\node (1) [label=above:{}] at (1,5) {};
\node (2) [label=above:{}] at (2,5) {};
\node (3) [label=above:{}] at (3,5) {};
\node (4) [label=above:{}] at (4,5) {};

\node (5) [label=above:{}] at (1,4) {};
\node (6) [label=above:{}] at (2,4) {};
\node (7) [label=above:{}] at (3,4) {};
\node (8) [label=above:{}] at (4,4) {};

\draw[very thick,blue] \convexpath{1,4,8,5}{.2cm};

\draw[blue,line width=1mm, ->] node[font=\fontsize{12}{0}\selectfont, right=of 4, left=-.5cm, below=0.2cm] {$~~~~~~~\rho_{b_1}$};

\node (9) [label=above:{}] at (1,3) {};
\node (10) [label=above:{}] at (2,3) {};
\node (11) [label=above:{}] at (3,3) {};
\node (12) [label=above:{}] at (4,3) {};

\node (13) [label=above:{}] at (1,2) {};
\node (14) [label=above:{}] at (2,2) {};
\node (15) [label=above:{}] at (3,2) {};
\node (16) [label=above:{}] at (4,2) {};

\node (17) [label=above:{}] at (1,1) {};
\node (18) [label=above:{}] at (2,1) {};
\node (19) [label=above:{}] at (3,1) {};
\node (20) [label=above:{}] at (4,1) {};

\draw[very thick,blue] \convexpath{9,12,20,17}{.2cm};

\draw[blue,line width=1mm, ->] node[font=\fontsize{12}{0}\selectfont, right=of 12, left=-.5cm, below=0.7cm] {$~~~~~~~\rho_{b_2}$};

  \begin{pgfonlayer}{background}
    \filldraw [line width=4mm,black!3]
      (0.2,0.2)  rectangle (4.8,5.8);
  \end{pgfonlayer}
\end{tikzpicture}}
\caption{Splitting of the partition $\rho$ with $\rho \vee \pi = \{\pi_1\cup\pi_2, \pi_3\cup\pi_4\cup\pi_5\}$ and $n=5$, $r=4$.}
\label{fig:diagram1-3}
\end{figure}

\begin{definition}
  Let $n , r \geq 1$,
  and let $\rho \in \Pi ([n]\times [r])$ be a partition of $[n]\times [r]$.
  \begin{enumerate}[\rm 1.] 
  \item
    For $b \subset [n]$, we let $\rho_b \subset \rho$ be defined as 
$$
\rho_b := \{ c \in \rho \ : \ c \subset b\times[r] \}. 
$$
\item Given $\eta \subset [n]$ 
 we split any partition $\rho$ of $\eta\times[r]$ into the equivalence classes
 deduced from the connected components of the graph $\rho_G$, as 
\begin{equation}
\label{fjklds1} 
\rho = \bigcup_{\substack{b\subset\eta\\b\times [r] \in \rho \vee \pi_\eta}} \rho_b. 
\end{equation} 
\end{enumerate}
\end{definition} 
 As an example, in Figure~\ref{fig:diagram1-3}-$a)$, when $b = \{1,2\}$ we have
$$
\rho_{\{1,2\}} = \big\{\{(1,1),(2,1),(2,2),(2,3)\},
\{(1,2),(1,3),(1,4),(2,4)\}\big\}, 
$$
and the partition \eqref{fjklds1} is
illustrated in Figure~\ref{fig:diagram1-3}-$b)$ with
 $b_1 = \{ 1,2\}$ and $b_2 = \{3,4,5\}$. 

\vspace{-0.4cm}
  
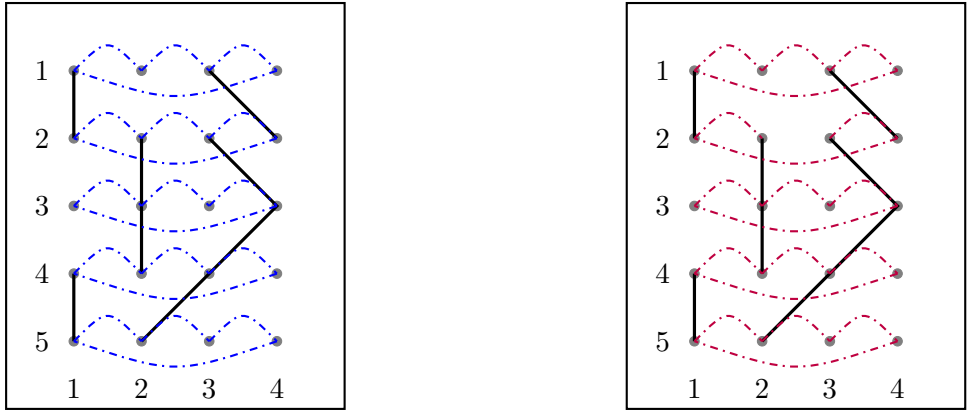
\begin{figure}[H]
\captionsetup[subfigure]{font=footnotesize}
\centering
\subcaptionbox{Diagram $\Gamma(\rho,\pi)$ and multigraph $\widetilde{\rho}_G$ in blue.}[.49\textwidth]{
\begin{tikzpicture}[scale=0.9] 
\draw[black, thick] (0,0) rectangle (5,6);

\node[anchor=east,font=\small] at (0.8,5) {1};
\node[anchor=east,font=\small] at (0.8,4) {2};
\node[anchor=east,font=\small] at (0.8,3) {3};
\node[anchor=east,font=\small] at (0.8,2) {4};
\node[anchor=east,font=\small] at (0.8,1) {5};

\node[anchor=south,font=\small] at (1,0) {1};
\node[anchor=south,font=\small] at (2,0) {2};
\node[anchor=south,font=\small] at (3,0) {3};
\node[anchor=south,font=\small] at (4,0) {4};

\filldraw [gray] (1,1) circle (2pt);
\filldraw [gray] (2,1) circle (2pt);
\filldraw [gray] (3,1) circle (2pt);
\filldraw [gray] (4,1) circle (2pt);
\filldraw [gray] (1,2) circle (2pt);
\filldraw [gray] (2,2) circle (2pt);
\filldraw [gray] (3,2) circle (2pt);
\filldraw [gray] (4,2) circle (2pt);
\filldraw [gray] (1,3) circle (2pt);
\filldraw [gray] (2,3) circle (2pt);
\filldraw [gray] (3,3) circle (2pt);
\filldraw [gray] (4,3) circle (2pt);
\filldraw [gray] (2,3) circle (2pt);
\filldraw [gray] (1,4) circle (2pt);
\filldraw [gray] (2,4) circle (2pt);
\filldraw [gray] (3,4) circle (2pt);
\filldraw [gray] (4,4) circle (2pt);
\filldraw [gray] (1,5) circle (2pt);
\filldraw [gray] (2,5) circle (2pt);
\filldraw [gray] (3,5) circle (2pt);
\filldraw [gray] (4,5) circle (2pt);

\draw[very thick] (1,5) -- (1,4); 
\draw[very thick] (3,5) -- (4,4);

\draw[very thick] (1,2) -- (1,1);
\draw[very thick] (2,2) -- (2,4);
\draw[very thick] (2,1) -- (3,2) -- (4,3) -- (3,4);

\foreach \i in {1,...,5}
         {
           \draw[thick,dash dot,blue] (1,\i) .. controls (1.5,\i+.5) .. (2,\i);
           \draw[thick,dash dot,blue] (2,\i) .. controls (2.5,\i+.5) .. (3,\i);
           \draw[thick,dash dot,blue] (3,\i) .. controls (3.5,\i+.5) .. (4,\i);
           \draw[thick,dash dot,blue] (4,\i) .. controls (2.5,\i-.5) .. (1,\i);
         }
         
  \begin{pgfonlayer}{background}
    \filldraw [line width=4mm,black!3]
      (0.2,0.2)  rectangle (4.8,5.8);
  \end{pgfonlayer}
\end{tikzpicture}}
\subcaptionbox{Diagram $\Gamma(\rho ,\pi)$ and graph $\rho_G$ in red.}[.49\textwidth]{
\begin{tikzpicture}[scale=0.9] 
\draw[black, thick] (0,0) rectangle (5,6);

\node[anchor=east,font=\small] at (0.8,5) {1};
\node[anchor=east,font=\small] at (0.8,4) {2};
\node[anchor=east,font=\small] at (0.8,3) {3};
\node[anchor=east,font=\small] at (0.8,2) {4};
\node[anchor=east,font=\small] at (0.8,1) {5};

\node[anchor=south,font=\small] at (1,0) {1};
\node[anchor=south,font=\small] at (2,0) {2};
\node[anchor=south,font=\small] at (3,0) {3};
\node[anchor=south,font=\small] at (4,0) {4};

\filldraw [gray] (1,1) circle (2pt);
\filldraw [gray] (2,1) circle (2pt);
\filldraw [gray] (3,1) circle (2pt);
\filldraw [gray] (4,1) circle (2pt);
\filldraw [gray] (1,2) circle (2pt);
\filldraw [gray] (2,2) circle (2pt);
\filldraw [gray] (3,2) circle (2pt);
\filldraw [gray] (4,2) circle (2pt);
\filldraw [gray] (1,3) circle (2pt);
\filldraw [gray] (2,3) circle (2pt);
\filldraw [gray] (3,3) circle (2pt);
\filldraw [gray] (4,3) circle (2pt);
\filldraw [gray] (2,3) circle (2pt);
\filldraw [gray] (1,4) circle (2pt);
\filldraw [gray] (2,4) circle (2pt);
\filldraw [gray] (3,4) circle (2pt);
\filldraw [gray] (4,4) circle (2pt);
\filldraw [gray] (1,5) circle (2pt);
\filldraw [gray] (2,5) circle (2pt);
\filldraw [gray] (3,5) circle (2pt);
\filldraw [gray] (4,5) circle (2pt);

\draw[very thick] (1,5) -- (1,4); 
\draw[very thick] (3,5) -- (4,4);

\draw[very thick] (1,2) -- (1,1);
\draw[very thick] (2,2) -- (2,4);
\draw[very thick] (2,1) -- (3,2) -- (4,3) -- (3,4);

\foreach \i in {1,...,3}
         {
           \draw[thick,dash dot,purple] (1,\i) .. controls (1.5,\i+.5) .. (2,\i);
           \draw[thick,dash dot,purple] (2,\i) .. controls (2.5,\i+.5) .. (3,\i);
           \draw[thick,dash dot,purple] (3,\i) .. controls (3.5,\i+.5) .. (4,\i);
           \draw[thick,dash dot,purple] (4,\i) .. controls (2.5,\i-.5) .. (1,\i);
         }

           \draw[thick,dash dot,purple] (1,4) .. controls (1.5,4+.5) .. (2,4);
           \draw[thick,dash dot,purple] (3,4) .. controls (3.5,4+.5) .. (4,4);
           \draw[thick,dash dot,purple] (4,4) .. controls (2.5,4-.5) .. (1,4);

\foreach \i in {5,...,5}
         {
           \draw[thick,dash dot,purple] (1,\i) .. controls (1.5,\i+.5) .. (2,\i);
           \draw[thick,dash dot,purple] (2,\i) .. controls (2.5,\i+.5) .. (3,\i);
           \draw[thick,dash dot,purple] (3,\i) .. controls (3.5,\i+.5) .. (4,\i);
           \draw[thick,dash dot,purple] (4,\i) .. controls (2.5,\i-.5) .. (1,\i);
         }
  \begin{pgfonlayer}{background}
    \filldraw [line width=4mm,black!3]
      (0.2,0.2)  rectangle (4.8,5.8);
  \end{pgfonlayer}
\end{tikzpicture}}
\caption{Connected non-flat partition diagram with $G$ a cycle graph and $n=5$, $r=4$.}
\label{fig:diagram2}
\end{figure}

\begin{definition}
 Let $n , r \geq 1$. 
 Given $\sigma \in \Pi ([n] )$ a partition of $[n]$, 
 we let $\Pi_\sigma ( [n] \times [r])$ denote the set of
 partitions $\rho$ of $[n] \times [r]$ such that
 $$
 \rho\vee\pi= \{ b \times [r] \ : \ b \in \sigma \}, 
 $$
 and
 we partition $\Pi ([n] \times [r])$ as
\begin{equation}
\label{partition} 
 \Pi ([n] \times [r])
 = \bigcup_{\sigma \in \Pi ([n] )} \Pi_\sigma ([n] \times [r]). 
\end{equation} 
\end{definition}
 We note that given $\eta \subset [n]$,
the set $\Pi_{\widehat{1}} ( \eta \times [r])$ consists of the partitions
$\rho$ of $\eta \times [r]$ for which the graph $\rho_G$
is connected, as in Figure~\ref{fig:diagram2}. 

\vskip-0.3cm

\noindent
The following lemma will be useful when applying induction
arguments on connected set partitions in $\Pi_{\widehat{1}}([n]\times [r])$.
\begin{lemma}
  \label{restrict-partition}
  Let $n\geq 2$.
  For any connected partition 
  $\rho\in \Pi_{\widehat{1}}( [n]\times [r])$
  there exists $i\in \{1,\dots,n \}$
  such that the set partition 
    $\{b\backslash\pi_i:b\in\rho\}$ 
  of $\{1,\dots,i-1,i+1,\dots,n \}\times [r]$
  is connected.
\end{lemma}
\begin{Proof}
 Let $\rho\in \Pi_{\widehat{1}}( [n]\times [r])$. 
 We consider the connected undirected graph
 $\eufrak{g}$ on $[n]$ in which two vertices $i,j\in [n]$
 are connected if and only
 if there exists a block $b\in \rho$ such that $\pi_i \cap b \not= \emptyset$
 and $\pi_j \cap b \not= \emptyset$,
 see Figure~\ref{fig:diagram0-11}-$a)$ for an example with $n=5$.  

\begin{figure}[H]
\captionsetup[subfigure]{font=footnotesize}
\centering
\subcaptionbox{Diagram $\Gamma(\rho,\pi)$ and graph $\eufrak{g}$.}[.48\textwidth]{
\begin{tikzpicture}[scale=0.9] 
\draw[black, thick] (0,0) rectangle (5,6);

\node[anchor=east,draw, circle, inner sep=0pt, minimum size=11pt,font=\small] at (0.7,5) {1};
\node[anchor=east,draw, circle, inner sep=0pt, minimum size=11pt,font=\small] at (0.7,4) {2};
\node[anchor=east,draw, circle, inner sep=0pt, minimum size=11pt,font=\small] at (0.7,3) {3};
\node[anchor=east,draw, circle, inner sep=0pt, minimum size=11pt,font=\small] at (0.7,2) {4};
\node[anchor=east,draw, circle, inner sep=0pt, minimum size=11pt,font=\small] at (0.7,1) {5};

\node[anchor=south,font=\small] at (1,0) {1};
\node[anchor=south,font=\small] at (2,0) {2};
\node[anchor=south,font=\small] at (3,0) {3};
\node[anchor=south,font=\small] at (4,0) {4};

\filldraw [gray] (1,1) circle (2pt);
\filldraw [gray] (2,1) circle (2pt);
\filldraw [gray] (3,1) circle (2pt);
\filldraw [gray] (4,1) circle (2pt);
\filldraw [gray] (1,2) circle (2pt);
\filldraw [gray] (2,2) circle (2pt);
\filldraw [gray] (3,2) circle (2pt);
\filldraw [gray] (4,2) circle (2pt);
\filldraw [gray] (1,3) circle (2pt);
\filldraw [gray] (2,3) circle (2pt);
\filldraw [gray] (3,3) circle (2pt);
\filldraw [gray] (4,3) circle (2pt);
\filldraw [gray] (2,3) circle (2pt);
\filldraw [gray] (1,4) circle (2pt);
\filldraw [gray] (2,4) circle (2pt);
\filldraw [gray] (3,4) circle (2pt);
\filldraw [gray] (4,4) circle (2pt);
\filldraw [gray] (1,5) circle (2pt);
\filldraw [gray] (2,5) circle (2pt);
\filldraw [gray] (3,5) circle (2pt);
\filldraw [gray] (4,5) circle (2pt);

\draw[very thick] (1,5) -- (1,4); 
\draw[very thick] (3,5) -- (4,4);

\draw[very thick] (1,2) -- (1,1);
\draw[very thick] (2,2) -- (2,4);
\draw[very thick] (2,1) -- (3,2) -- (4,3) -- (3,4);

\draw[very thick,purple] (0.3,5) .. controls (0.2,4.5) .. (0.3,4);
\draw[very thick,purple] (0.3,4) .. controls (0.2,3.5) .. (0.3,3);
\draw[very thick,purple] (0.3,3) .. controls (0.2,2.5) .. (0.3,2);
\draw[very thick,purple] (0.3,2) .. controls (0.2,1.5) .. (0.3,1);
\draw[very thick,purple] (0.7,4) .. controls (0.85,3) .. (0.7,2);
\draw[very thick,purple] (0.7,4) .. controls (0.95,3) .. (0.7,1);

  \begin{pgfonlayer}{background}
    \filldraw [line width=4mm,black!3]
      (0.2,0.2)  rectangle (4.8,5.8);
  \end{pgfonlayer}
\end{tikzpicture}
}
\subcaptionbox{Diagram $\Gamma(\rho,\pi)$ and spanning tree $\widebar{\eufrak{g}}$.}[.49\textwidth]{
\begin{tikzpicture}[scale=0.9] 
\draw[black, thick] (0,0) rectangle (5,6);

\node[anchor=east,draw, circle, inner sep=0pt, minimum size=11pt, font=\small] at (0.7,5) {1};
\node[anchor=east,draw, circle, inner sep=0pt, minimum size=11pt,font=\small] at (0.7,4) {2};
\node[anchor=east,draw, circle, inner sep=0pt, minimum size=11pt,font=\small] at (0.7,3) {3};
\node[anchor=east,draw, circle, inner sep=0pt, minimum size=11pt,font=\small] at (0.7,2) {4};
\node[anchor=east,draw, circle, inner sep=0pt, minimum size=11pt,font=\small] at (0.7,1) {5};

\node[anchor=south,font=\small] at (1,0) {1};
\node[anchor=south,font=\small] at (2,0) {2};
\node[anchor=south,font=\small] at (3,0) {3};
\node[anchor=south,font=\small] at (4,0) {4};

\filldraw [gray] (1,1) circle (2pt);
\filldraw [gray] (2,1) circle (2pt);
\filldraw [gray] (3,1) circle (2pt);
\filldraw [gray] (4,1) circle (2pt);

\filldraw [gray] (1,2) circle (2pt);
\filldraw [gray] (2,2) circle (2pt);
\filldraw [gray] (3,2) circle (2pt);
\filldraw [gray] (4,2) circle (2pt);

\filldraw [gray] (1,3) circle (2pt);

\filldraw [gray] (2,3) circle (2pt);
\filldraw [gray] (3,3) circle (2pt);

\filldraw [gray] (4,3) circle (2pt);
\filldraw [gray] (2,3) circle (2pt);
\filldraw [gray] (1,4) circle (2pt);
\filldraw [gray] (2,4) circle (2pt);
\filldraw [gray] (3,4) circle (2pt);
\filldraw [gray] (4,4) circle (2pt);
\filldraw [gray] (1,5) circle (2pt);
\filldraw [gray] (2,5) circle (2pt);

\filldraw [gray] (3,5) circle (2pt);
\filldraw [gray] (4,5) circle (2pt);

\draw[very thick] (1,5) -- (1,4); 
\draw[very thick] (3,5) -- (4,4);

\draw[very thick] (1,2) -- (1,1);
\draw[very thick] (2,2) -- (2,4);
\draw[very thick] (2,1) -- (3,2) -- (4,3) -- (3,4);

\draw[very thick,blue] (0.3,5) .. controls (0.2,4.5) .. (0.3,4);

\draw[very thick,blue] (0.3,3) .. controls (0.2,2.5) .. (0.3,2);
\draw[very thick,blue] (0.3,2) .. controls (0.2,1.5) .. (0.3,1);
\draw[very thick,blue] (0.7,4) .. controls (0.85,3) .. (0.7,2);

  \begin{pgfonlayer}{background}
    \filldraw [line width=4mm,black!3]
      (0.2,0.2)  rectangle (4.8,5.8);
  \end{pgfonlayer}
\end{tikzpicture}
}
\caption{Example of graph $\eufrak{g}$ and its spanning tree subgraph.}
\label{fig:diagram0-11}
\end{figure}
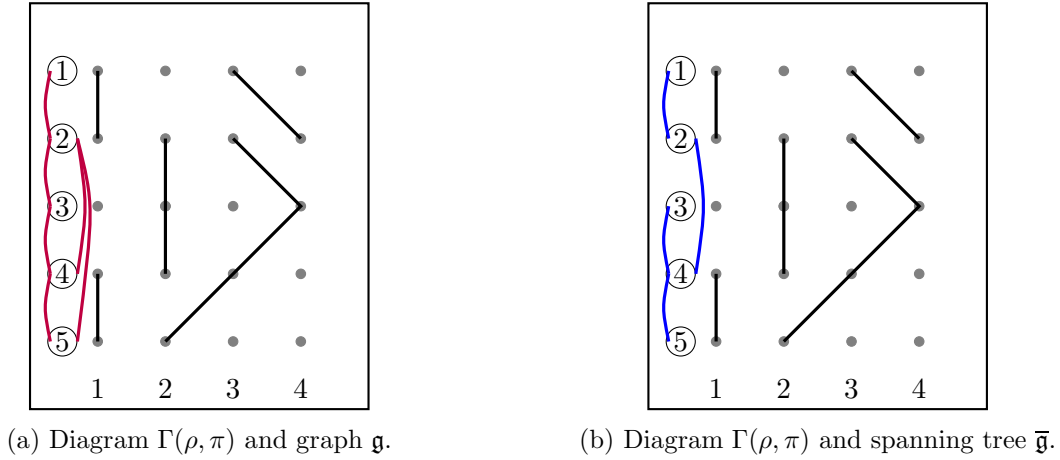

\vspace{-.4cm}
  
\noindent
By e.g. Theorem~4.2.3 in \cite{balakrishnan},
 $\eufrak{g}$ contains a spanning
 tree $\widebar{\eufrak{g}}$, 
 as shown in Figure~\ref{fig:diagram0-11}-$b)$. 
 Let $i\in [n]$ be a leaf in the tree $\widebar{\eufrak{g}}$. 
 If the partition 
 ${\rho}^{(i)}:=\{ b \setminus \pi_i \ \!  : \ \!  b\in \rho \}$
 of $([n] \setminus \{i\} ) \times [r]$ 
 had more than one connected component, then, for $\rho$ to be connected,
 $\pi_i$ would have to connect to all such components, 
 hence the vertex $i$ would be adjacent to
 more than one vertex in $\widebar{\eufrak{g}}$,
 which is not the case.
\end{Proof}
 In what follows, we say that a partition $\rho$ of $[n]\times [r]$ 
 is {\em non-flat} if 
 its diagram $\Gamma(\rho,\pi )$ is non-flat, i.e. 
 if $\rho \wedge \pi = \widehat{0}$,
 see Chapter~4 of \cite{peccatitaqqu} and Figure~\ref{fig:diagram2}. 

\begin{definition} 
  \noindent
 Given $n,r\geq 1$, we denote by 
$$
 {\rm NF} (n,r) :=\{\rho\in\Pi ([n]\times[r]) \ : \ \rho\wedge \pi=\widehat{0} \}
  $$
 the set of 
 non-flat partitions of $[n]\times[r]$, and by
$$
 {\rm CNF} (n,r) :=\{\rho\in\Pi_{\widehat{1}} ([n]\times[r]) \ : \ \rho\wedge \pi=\widehat{0} \}
  $$
 the set of connected
 non-flat partitions of $[n]\times[r]$. 
\end{definition}
\noindent 
We will also consider the following set of 
 connected non-flat partitions which have a maximal
 number of blocks. 
 \begin{definition} 
  \noindent
 Given $n\geq 1$ and $r\geq 2$,
we denote by
$$
 {\rm M}(n,r):=\{\rho\in {\rm CNF} (n,r)  \ : \ |\rho|= 1 + (r-1)n \}
$$
 the set of maximal connected non-flat partitions of $[n]\times[r]$.
\end{definition}

\noindent
 The bound in part $(a)$ of the next lemma is consistent with
 (6.2) in Proposition~6.1 of \cite{schulte-thaele},
 which shows that the power $r$ of $n!$ cannot be improved in \eqref{coeff-0}. 
\begin{lemma}
  \label{fjkldsf-l}
  \noindent
  $a)$ The cardinality of the set \ ${\rm NF} (n,r)$
  of 
  non-flat partitions of $[n]\times[r]$ satisfies 
 \begin{equation}
   \label{coeff-0}
  |{\rm NF} (n,r) | \leq n!^r r!^{n-1}, 
  \qquad n,r \geq 1. 
\end{equation}
\noindent
$b)$ 
 The cardinality of the set 
${\rm M}(n,r)$
 of maximal connected non-flat partitions of $[n]\times[r]$ satisfies 
\begin{equation}\label{coeff-10}
  |{\rm M}(n,r)|=r^{n-1}\prod_{i=1}^{n-1}(1+(r-1)i),
  \qquad n,r\geq 1, 
\end{equation}
 with the bounds 
\begin{equation}\label{coeff-1}
    ( (r-1)r )^{n-1}(n-1)!\le
    |{\rm M}(n,r)|
     \leq ( (r-1)r )^{n-1}n!, \quad n\geq 1, \ r\geq 2. 
\end{equation}
\end{lemma}
\begin{Proof}
\noindent
 $a)$ We clearly have $|{\rm NF} (1,r) | =1$ for all $r\geq 1$. 
Any 
 non-flat partition $\rho \in {\rm NF}(n+1,r)$ can be
 obtained from a 
  non-flat partition in ${\rm NF}(n,r)$ in at most
 $(n+1)^r r!$ ways, by connecting each of the $r$ new points
 in at most $n+1$ possible ways (including non-connection),
 and multiplying by the number $r!$ of possible permutations of $\pi_i$. 
 This yields the induction inequality 
 $$
 |{\rm NF} (n+1,r) | \leq (n+1)^r r! | {\rm NF} (n,r) |,
 $$
 from which we conclude that \eqref{coeff-0} holds. 
 
 \medskip 

 \noindent
 $b)$ 
 We clearly have $|{\rm M}(1,r)|=1$ for all $r\geq 1$. 
 Next, each maximal connected non-flat partition 
 $\rho\in {\rm M}(n+1,r)$ can be
   obtained by choosing 
   one of $r$ elements of $\{(n+1,1),\dots,(n+1,r)\}$, 
   and connecting them in $1+(r-1)n$ ways 
   to any partition in ${\rm M}(n,r)$, $n\geq 1$.
   This implies the recursion formula 
   $$|{\rm M}(n+1,r)|=r\times (1+(r-1)n)|{\rm M}(n,r)|,
   $$
 which yields \eqref{coeff-10}.
\end{Proof} 

\section{Virtual cumulants} 
\label{s3}
\noindent
The following definition
uses the concept of independence of a virtual field
with respect to graph connectedness,
see Relation~(17) in \cite[p.~34]{MalyshevMinlos91}. 
\begin{definition}
  Let $n,r\geq 1$.
  We say that a mapping $F$ defined on partitions of
  $[n]\times [r]$
    admits the {\it connectedness factorization} property
  if it decomposes according to the partition \eqref{fjklds1} as 
  \begin{equation}
    \label{dia-factoriz}
  F ( \rho ) = \prod_{b \times [r] \in \rho \vee \pi } F ( \rho_b ),  
  \qquad
   \rho \in \Pi ([n]\times [r]). 
\end{equation}
\end{definition} 
In what follows,
given $F$ a mapping defined on the partitions of
  $[n]\times [r]$, 
  we will use the M\"obius transform
  $\widehat{F}$ of $F$, 
  defined as
  \begin{equation}
\nonumber 
  \widehat{F}( \eta ):=\sum_{\rho \in \Pi ( \eta \times [r])} F(\rho ),
  \qquad
  \eta \subset [n],
\end{equation}
with $\widehat{F}(\emptyset):=0$,
see \cite{rota1964} 
 and \S~2.5 of \cite{peccatitaqqu}. 
  We refer to \cite[p.~33]{MalyshevMinlos91} for the following definition.
\begin{definition}
  Let $n,r\geq 1$.
  The virtual cumulant $G$ of 
 a mapping $F$ on $\bigcup_{\eta \subset [n]} \Pi ( \eta \times [r])$ 
 is defined by letting  
 $C_F ( \eta ):=\widehat{F}( \eta )$ when $| \eta |=1$, and
 then recursively by 
\begin{equation}
  \label{fjkldsf}
  C_F ( \eta ):=\widehat{F}( \eta )-\sum_{
    { \sigma \in \Pi (\eta )
    \atop 
    |\sigma|\geq 2
    }
  }
  \prod_{b \in \sigma } C_F ( b ), \qquad
  \eta \subset [n], \quad |\eta |\geq 2. 
\end{equation}
\end{definition} 
Relation~\eqref{fjkldsf} also implies the relation 
 \begin{equation}
   \label{jlkdf3} 
 C_F ( \eta )=
 \sum_{
   \sigma \in \Pi ( \eta )
 }
 (-1)^{|\sigma |-1} (|\sigma |-1)! \prod_{b\in \sigma} \widehat{F} ( b ), 
\end{equation}
 see Relation~(16') page~33 of \cite{MalyshevMinlos91},
 which is also the classical cumulant-moment relationship,
 see e.g. Corollary~3.2.2 in \cite{peccatitaqqu}. 
 The following proposition is an extension of
the classical Lemma~2 in \cite[p.~34]{MalyshevMinlos91},
see also Lemma~3.1 in \cite{khorunzhiy}.
\begin{prop}
\label{mainthm-1}
Let $n,r\geq 1$. 
Let $F$ be a mapping defined on $\bigcup_{\eta \subset [n]}
 \Pi ( \eta \times [r])$ and admitting the {connectedness factorization} property
 \eqref{dia-factoriz}.
 Then, for $\eta \subset [n]$ with $\eta \ne\emptyset$, 
 the virtual cumulant of $F$ is given by the sum 
\begin{equation}
\label{keylma-2}
C_F (\eta ) = \sum_{\sigma \in \Pi_{\widehat{1}} ( \eta \times [r] )
  \atop {\rm (connected)}}
F(\sigma )
\end{equation}
over connected partitions on $\eta \times [r]$. 
\end{prop} 
\begin{Proof}
The claim is true when $|\eta |=1$. 
Assume that it is true for all $\eta \subset [n]$ for some $n \geq 1$, and 
 let $\eta$ be such that $|\eta |=n+1$.
 By \eqref{partition} and \eqref{dia-factoriz}, we have 
\begin{eqnarray*}
  \widehat{F}(\eta )
  &=&
  \sum_{\rho \in \Pi (\eta \times [r])} F(\rho )
  \\
  &=&
  \sum_{\sigma \in \Pi (\eta )}
  \sum_{\rho \in \Pi_\sigma (\eta \times [r])} F(\rho )
  \\
  &=&
  \sum_{\sigma \in \Pi (\eta )}
  \sum_{\rho \in \Pi_\sigma (\eta \times [r])}
  \prod_{b \in \sigma} 
  F(\rho_b ) 
 \\
  &=&
  \sum_{\sigma \in \Pi (\eta )}
  \prod_{b \in \sigma } 
  \sum_{\rho \in \Pi_{\widehat{1}} (b \times [r])
    \atop {\rm (connected)}
  } 
  F(\rho )
 \\
  &=&
 \sum_{\rho \in \Pi_{\widehat{1}} ( \eta \times [r] )
   \atop {\rm (connected)}
 }
  F(\rho )
  +
  \sum_{\sigma \in \Pi (\eta ) \atop |\sigma |\geq 2}
  \prod_{b \in \sigma } 
  C_F (b ), 
\end{eqnarray*}
where the last equality follows from the induction hypothesis \eqref{keylma-2} 
when $| \eta |\leq n$.
 The proof is completed by subtracting the last term on both sides. 
\end{Proof}
\section{Cumulants of multiparameter stochastic integrals} 
\label{s4}
\noindent 
 Proposition~\ref{p01-1-4} rewrites the product in \eqref{nthexpectation0-0} 
 of Proposition~\ref{p01-1-0} 
 as a product on the edges of the graph $\rho_G$ 
 similarly to Proposition~4 in \cite{prkhp} when $v(x,y)$ vanishes on
 the diagonal, and it generalizes Proposition~2.4 of \cite{jansen}
 from two-parameter Poisson stochastic
 integrals to multiparameter integrals of higher orders. 
\begin{prop}
\label{p01-1-4} 
Let $n \geq 1$, $r\geq 2$, and assume that the process $v(x,y)$ vanishes
on diagonals, i.e. $v(x,x) = 0$, $x\in \real^d$. 
Then, the $n$-$th$ moment of the multiparameter stochastic integral
 \eqref{e1} is given by the summation 
\begin{equation} 
\nonumber 
\sum_{\rho \in \Pi ( [n] \times [r] )
\atop {
    \rho \wedge \pi = \widehat{0}
    \atop
        {\rm (non-flat)}
}
}
\int_{(\R^d )^{|\rho |}}
  \prod_{\{\eta_1,\eta_2\} \in E(\rho_G)}
 \E \big[ v(x_{\eta_1},x_{\eta_2})^{m({\eta_1,\eta_2})} \big] 
   \ \! \prod_{\eta \in V(\rho_G)}
\Lambda(\mathrm{d}x_\eta ), 
\end{equation} 
over connected non-flat partitions,
where $m ({\eta_1,\eta_2})$ represents the multiplicity of the edge
$(\eta_1,\eta_2)$ in the multigraph $\widetilde{\rho}_G$. 
\end{prop}
The next proposition is a consequence of Propositions~\ref{mainthm-1}
and \ref{p01-1-4},
 and it also extends Proposition~2.5 of \cite{jansen}
 from the two-parameter case to the multiparameter
 case. 
 Note that in our setting, the two-parameter case
 only applies to the edge counting. 
\begin{prop}
\label{p01-1} 
Let $n \geq 1$, $r \geq 2$, and assume that the process $v(x,y)$ vanishes
on diagonals, i.e. $v(x,x) = 0$, $x\in \real^d$. 
Then, the $n$-$th$ cumulant of the multiparameter stochastic integral
 \eqref{e1} is given by the summation 
\begin{equation}
  \label{nthexpectation0-2}
    \sum_{\rho \in \Pi_{\widehat{1}} ( [n] \times [r])
    \atop
    {\rho \wedge \pi = \widehat{0}
    \atop {\rm (non-flat \ \! connected)}
    }
    }
    \int_{(\R^d )^{|\rho |}}
      \prod_{\{\eta_1,\eta_2\} \in E(\rho_G)}
  \E \big[ v (x_{\eta_1},x_{\eta_2})^{m({\eta_1,\eta_2})} \big] 
 \ \! \prod_{\eta \in V(\rho_G)}
\Lambda(\mathrm{d}x_\eta )
\end{equation}
over connected non-flat partitions.
\end{prop}
\begin{Proof}
 The functional 
\begin{equation}
\nonumber 
F(\rho ):=
\int_{(\R^d )^{|\rho |}}
  \prod_{\{\eta_1,\eta_2\} \in E(\rho_G)}
 \E \big[ v(x_{\eta_1},x_{\eta_2})^{m({\eta_1,\eta_2})} \big] 
   \ \! \prod_{\eta \in V(\rho_G)}
\Lambda(\mathrm{d}x_\eta )
\end{equation}
 satisfies the connectedness factorization property
 \eqref{dia-factoriz}, as for
 $\sigma = b \times [r] \in \rho \vee \pi$ and
 $\sigma' = b' \times [r] \in \rho \vee \pi$
 with $b \not= b'$, the variables $(x_\eta )_{\eta \in \rho_b}$
 are distinct from the variables $(x_\eta )_{\eta \in \rho_{b'}}$
 in the above integration. 
 Hence, Relation~\eqref{nthexpectation0-2} follows from Proposition~\ref{mainthm-1}
 and the classical cumulant-moment relationship \eqref{jlkdf3},
 since by Proposition~\ref{p01-1-0},
 $\widehat{F}([n])$ is the $n$-$th$ moment of the multiparameter stochastic integral \eqref{e1}. 
\end{Proof}
\section{Cumulants of subgraph counts} 
\label{s5}
\noindent
 Let 
 $H:\R^d\times\R^d\to[0,1]$ 
 denote a 
 measurable connection function such that
 \begin{equation}
 \nonumber 
0<\int_{\real^d}H(x,y) \Lambda(\mathrm{d}x)< \infty, 
\end{equation}
 for all $y \in\R$. 
 Given $\omega \in \Omega$, for any $x,y\in\omega$ with $x\not= y$,
 an edge connecting $x$ and $y$ is added with probability $H(x,y)$,
 independently of the other pairs.
  The resulting random graph, together with the
 {Poisson point process $\Xi$
 with intensity measure $\Lambda$ on $\real^d$}, 
 is called the random-connection model and denoted by $G_H(\Xi)$.
 
 \medskip

 In the case where the connection function $H$ is given by $H(x,y) := \bone_{\{\|x-y\|\leq R\}}$ for some $R>0$, the resulting graph is completely determined by the geometry of the underlying point process $\Xi$, and is called a random geometric graph {\cite{penrosebk}}, 
 which is included as a special case in this paper. The main case of interest in this section is the general random-connection model. We shall have further discussion about the random geometric graph in Section~\ref{rgg}, as we believe it is of independent interest. 

\medskip

 Given $G$ a connected graph with $|V(G)| = r$ vertices,
 we denote $N_G$ the count of subgraphs isomorphic to $G$ in the
 random-connection model $G_H( \Xi )$,
 which can be represented as the multiparameter stochastic integral 
\begin{equation}
\nonumber 
  N_G:=\sum_{\{ V_1,\dots,V_r\} \subset \omega  } \ 
  \prod_{\{i,j\} \in E(G)} 
  \bone_{\{V_i\leftrightarrow V_j\}}
  = \int_{(\real^d)^r}
  \prod_{\{i,j\} \in E(G)} 
  \bone_{\{x_i\leftrightarrow x_j\}}
  \ \! \omega (\mathrm{d}x_1)\cdots \omega (\mathrm{d}x_r), 
\end{equation}
up to division by the number of automorphisms of $G$.
Here, 
$\bone_{\{x \leftrightarrow y\}}$ denotes a $\{0,1\}$-valued
Bernoulli random variable with parameter $H(x,y)$ 
 when $x\not=y$, $x,y\in \real^d$, and 
\begin{equation}
\label{jflds} 
\bone_{\{x\leftrightarrow x\}}:=0, \qquad x\in \real^d. 
\end{equation}
 Consequently, we have $\bone_{\{V_i\leftrightarrow V_j\}}=1$ or $0$ depending
 on whether $V_i$ and $V_j$ are connected or not by an edge
in $G_H(\Xi)$. 
\noindent
 The first moment of $N_G$ can be computed as 
\begin{equation}
\label{fm} 
\E [ N_G ] =\int_{(\R^d)^r}\left(\prod_{\{i,j\}\in E(G)}H(x_i,x_j)\right)\prod_{i=1}^r\Lambda(\mathrm{d}x_i).
\end{equation}
 Higher order moments can be computed from the following result which
is a direct consequence of Proposition~\ref{p01-1}
by taking $v(x,y) : = \bone_{\{x \leftrightarrow y\}}$
in \eqref{nthexpectation0-2} and by using
{\em non-flat} partition diagrams $\Gamma(\rho,\pi )$
 such that $\rho \wedge \pi = \widehat{0}$, 
 to take into account condition~\eqref{jflds}.  
\begin{prop}
 \label{lma-diagram1}
  Let $n \geq 1$ and $r\geq 2$.
  The moments and cumulants of $N_G$ are given by the summation 
\begin{equation} 
\label{nthexpectation-2} 
\E [(N_G)^n] =
\sum_{\rho \in \Pi ( [n] \times [r] )
  \atop {
    \rho \wedge \pi = \widehat{0}
    \atop
        {\rm (non-flat)}
  }
}
\int_{(\R^d )^{|\rho |}}
\Bigg(
\prod_{\{\eta_1,\eta_2\} \in E(\rho_G)}
H(x_{\eta_1},x_{\eta_2})
\Bigg)
\prod_{\eta \in V(\rho_G)}
\Lambda (\mathrm{d}x_\eta ), 
\end{equation} 
 over non-flat partitions, and by the summation   
\begin{equation}
  \label{cumulant-diagram1}
  \kappa_n(N_G)=\sum_{\rho \in \Pi_{\widehat{1}} ( [n] \times [r])
    \atop
    {\rho \wedge \pi = \widehat{0}
    \atop {\rm (non-flat \ \! connected)}}
  }
  \int_{(\R^d )^{|\rho |}}
  \Bigg(
\prod_{\{\eta_1,\eta_2\} \in E(\rho_G)}
H(x_{\eta_1},x_{\eta_2})
\Bigg)
\prod_{\eta \in V(\rho_G)}
\Lambda (\mathrm{d}x_\eta ), 
\end{equation} 
 over connected non-flat partitions. 
\end{prop}
\begin{Proof}
 Relations~\eqref{nthexpectation-2}-\eqref{cumulant-diagram1}
 are consequence of Proposition~\ref{p01-1}, after taking
  $v (x_i,x_j):=\bone_{\{x_i\leftrightarrow x_j\}}$,
  $\{i,j\}\in E(G)$. 
  The summations are restricted to {\em non-flat} partitions
  due to condition~\eqref{jflds} as in Section~2 of \cite{prkhp}. 
\end{Proof}
\section{Asymptotic growth of subgraph count cumulants}
\label{s6}
\noindent 
In this section we consider the following assumption,
where $(\Lambda_\lambda)_{\lambda >0}$ is a family 
of {$\sigma$-finite} intensity measures on $\real^d$ and $H(x,y)$ is the connection function of the random-connection model.
\begin{assumption}
\label{a61} 
 Let $r\geq 2$ and $n \geq 1$.
 There exist constants $c_H, C_H > 0$
 such that for any connected non-flat partition 
 $\rho \in \Pi_{\widehat{1}} ( [n] \times [r])$, we have 
\begin{equation}
\label{integ-connecting3}
  c_H^{|E(\rho_G)|} ( \lambda C_H )^{|V(\rho_G)|}
\leq     \int_{\R^d}\cdots\int_{\R^d}
\Bigg(
\prod_{\{i,j\}\in E(\rho_G) }H(x_i,x_j)
\Bigg)
\
\prod_{k\in V(\rho_G) } \Lambda_\lambda (\mathrm{d}x_k),
\quad \lambda >0. 
\end{equation}
\end{assumption}
 We consider two settings satisfying Assumption~\ref{a61}.

 \vspace{-0.2cm}

 \noindent
\begin{example}
[Increasing intensity] When the intensity measure $\Lambda_\lambda$ takes the form 
$$\Lambda_\lambda (\mathrm{d}x)=\lambda\mu(\mathrm{d}x),\qquad\lambda>0,$$
 for $\mu$ a finite diffuse measure on $\R^d$, 
 Assumption~\ref{a61} is satisfied by any translation-invariant
 continuous kernel function $H : \real^d\times \real^d \to [0,1]$ non vanishing at $(0,0)$. 
 Indeed, in this case there exist $c_H>0$ and a Borel set $B\subset \real^d$
 such that $\mu ( B)>0$ and
$$
H(x,y) = H(x-y,0) \geq c_H \bone_B(x)\bone_B(y), \qquad x,y \in \real^d,
$$
hence 
\begin{eqnarray*} 
  c_H^{|E(\rho_G)|} ( \mu ( B ))^{|V(\rho_G)|} & = & 
  c_H^{|E(\rho_G)|}
  \int_B \cdots \int_B 
  \prod_{k\in V(\rho_G) } \mu (\mathrm{d}x_k)
  \\
   & \leq &  
  \int_{\R^d}\cdots\int_{\R^d}
  \Bigg(
  \prod_{\{i,j\}\in E(\rho_G)} H(x_i,x_j)
  \Bigg)
  \prod_{k\in V(\rho_G) } \mu (\mathrm{d}x_k), 
\end{eqnarray*}
so that we can take $C_H := \mu (B)$. 
\end{example}
 The setting of the {above increasing intensity example 
 covers the following long and short range dependence settings}: 
\begin{enumerate}[i)]
\item
 the power-law fading kernel 
 $H(x,y) = 1 \wedge \Vert x - y \Vert^{- \beta}$, $x,y\in \real^d$, for some $\beta > 0$,  
\item
 the Rayleigh fading kernel $H(x,y) = e^{ - \beta \Vert x - y\Vert^2}$, $x,y\in \real^d$, for some $\beta > 0$, 
\item
 the Boolean kernel $H(x,y)=\bone_{\{\|x-y\|\leq R \}}$, $x,y\in \real^d$,
 for some $R >0$, which yields the random geometric graph, with
 $C_H = v_d (R /2)^d$,
 where $v_d$ denotes the volume of the unit ball in $\real^d$. 
\end{enumerate}
\vspace{-0.2cm}
\noindent
\begin{example}
[Growing observation window] When the intensity measure $\Lambda_\lambda$ takes the form 
$$\Lambda_\lambda (\mathrm{d}x)= {\bf 1}_{A_\lambda} (x) \mu(\mathrm{d}x),\qquad\lambda>0,$$
  where $(A_\lambda )_{\lambda >0}$ is a non-decreasing sequence
  of Borel subsets of $\real^d$
  {such that $A_\lambda \uparrow \real^d$}
  as $\lambda$ tends to infinity, 
  $\mu$ a 
  diffuse
  {$\sigma$-finite}
  measure on $\R^d$, 
  and the kernel function $H : \real^d\times \real^d \to [0,1]$
  is lower bounded as $H(x,y) \geq c_H$, $x,y\in \real^d$ for some $c_H>0$,
  Assumption~\ref{a61} is satisfied with 
  $\mu ( A_\lambda ) = C_H \lambda$,
  e.g. when $A_\lambda$ is  the ball of radius $\lambda^{1/d} >0$
  and $\mu$ is the Lebesgue measure on $\real^d$,
  with $C_H = v_d$. 
\end{example} 

\noindent
 Next, we investigate the asymptotic behaviour of the cumulants $\kappa_n(N_G)$
as the intensity $\lambda$ tends to infinity, as a consequence of the partition diagram representation of cumulants. 
In what follows, given two positive functions $f$ and $g$ on $(1,\infty )$ we write
$f(\lambda ) \ll g(\lambda )$ if 
$\lim_{\lambda \to \infty} g(\lambda ) / f(\lambda ) = \infty$.
\begin{definition}
  Let $G$ be a connected graph with $|V(G)| = r$ vertices, $r\geq 2$.
  For every $\lambda >0$,
  let $G_{H_\lambda} (\Xi)$
  denote the random-connection model
  with connection function 
  $$
  H_\lambda(x,y):= c_\lambda H(x,y), \qquad x,y \in \real^d.
  $$ 
  We consider the following regimes.
\begin{itemize}
\item Dilute regime: for some constant $K>0$ we have 
\begin{equation}
    \label{fjnldsf}
    \lambda^{-1/\zeta} \ll c_\lambda \leq K,
    \qquad
    \lambda\to \infty,
\end{equation}
 where
 \begin{equation}
   \label{fjkldf}
   \zeta:=\max\left\{\frac{|E(H)|}{|V(H)|-1} \ \! : \ \!
   H\subseteq G, \ \! |V(H)| \geq 2 \right\}.
   \end{equation} 
\item Sparse regime: for some constants $K>0$ and $\alpha \geq 1$ we have
  \begin{equation}
    \label{fjnldsf-2}
    c_\lambda \leq \frac{K}{\lambda^\alpha},
        \qquad
    \lambda\to \infty. 
    \end{equation} 
\end{itemize} 
\end{definition}
In case $c_\lambda = K$ for all $\lambda > 0$
we also say that we are in the full random graph regime,
and in the sequel we take $K=1$ for simplicity.
We note that in general we have $\zeta > 1$ except when $G$ is a tree,
in which case $\zeta = 1$. 
\begin{prop}[Dilute regime]
\label{t1}
  Let $G$ be a connected graph with $|V(G)| = r$ vertices, $r\geq 2$, 
satisfying
Assumption~\ref{a61} for all $n\geq 1$
in the dilute regime \eqref{fjnldsf}. 
  We have the cumulant bounds 
\begin{equation}
\label{equiv-1}
 (n-1)! c_\lambda^{n |E(G)| } ( K_1 \lambda )^{1+(r-1)n}
 \leq 
 \kappa_n(N_G)
 \leq 
 n!^r c_\lambda^{n |E(G) |} ( K_2 \lambda )^{1+(r-1)n},
 \quad \lambda \geq 1, 
\end{equation}
for some constants $K_1$, $K_2>0$ independent of $\lambda$ and $n\geq 1$.
\end{prop} 
\begin{Proof}
   We identify the leading terms in the sum \eqref{cumulant-diagram1}
 over connected non-flat partitions, in which 
 every summand involves 
 a factor $c_\lambda^{|E(\rho_G)|} \lambda^{|V(\rho_G)|}$,
 since every
 vertex in $\rho_G$ contributes a factor $\lambda$, and that 
 every edge contributes a factor $c_\lambda$.
 Therefore, it is sufficient to show that for any $\rho_G$ 
 with $\rho \in {\rm CNF} (n,r)$
 a connected non-flat partition of $[n]\times[r]$, 
   we have
  \begin{equation}\label{negli}
    c_\lambda^{|E(\rho_G)|} \lambda^{|V(\rho_G)|}
    = O \big(\lambda^{1+(r-1)n}c_\lambda^{n |E(G)|}\big). 
\end{equation}
 This claim follows from \eqref{fm} when $n=1$.
Suppose that \eqref{negli} holds up to the rank $n\ge1$, 
 and let $\rho \in \Pi_{\widehat{1}}( [n+1]\times [r])$ 
 be a connected non-flat partition. 
By Lemma~\ref{restrict-partition}, there exists $i\in[n+1]$ such that
  the subgraph $\widebar{\rho}_G$
  induced by $\rho_G$ on the vertex set
$$
V(\widebar{\rho}_G):=
\big\{
b\in \rho \ \!  : \ \!  b \cap (\cup_{j\ne i}\pi_j) \not= \emptyset
\big\} 
$$
is connected. 
 Let $\widehat{\rho}_G $ 
 denote the subgraph induced by $\rho_G$ on the vertex set
$$
V(\widehat{\rho}_G):=
\big\{
b\in \rho \ \!  : \ \!  b \cap \pi_{i} \not= \emptyset
\big\}, 
$$
with
 $\widehat{\rho}_G\simeq G$
because $\rho$ is non-flat, and let $H:=\widebar{\rho}_G\cap\widehat{\rho}_G$. 
Since $H\subseteq\widehat{\rho}_G$ 
we have  
$$
\lambda^{|V(H)|-1} c_{\lambda}^{|E(H)|} \geq 
\big( \lambda c_{\lambda}^{\zeta }\big)^{|E(H)|/\zeta },
\qquad
\lambda \geq 1. 
$$ 
 Hence from 
$
  \lim_{\lambda \to \infty} 
  \lambda c_\lambda^{\zeta} = \infty$ 
  we get 
  $$
  \liminf_{\lambda \to \infty}
  \lambda^{|V(H)|-1} c_{\lambda}^{|E(H)|} 
  >0. 
  $$ 
 On the other hand, by the induction hypothesis we have 
\begin{equation}
  \frac{\lambda^{|V(\widebar{\rho}_G)|}c_\lambda^{|E(\widebar{\rho}_G)|}}{\lambda^{1 + (r-1)n }c_\lambda^{n|E(G)|}} = O(1), 
\end{equation}
hence, since
$|V(\widehat{\rho}_G)| = |V(G)|$
and
$|E(\widehat{\rho}_G)| = |E(G)|$, 
\begin{eqnarray}
  \nonumber
  \frac{\lambda^{|\rho|}c_\lambda^{|E(\rho_G)|}}{\lambda^{1 + (r-1)(n+1) }c_\lambda^{(n+1)|E(G)|}}
  &=&\frac{\lambda^{|V(\widebar{\rho}_G)|+|V(\widehat{\rho}_G)|-|V(H)|}c_\lambda^{|E(\widebar{\rho}_G)|+|E(\widehat{\rho}_G)|-|E(H)|}}{\lambda^{1 + (r-1)(n+1) }c_\lambda^{(n+1)|E(G)|}}
  \\
  \nonumber
  &=&\frac{\lambda^{|V(\widebar{\rho}_G)|}c_\lambda^{|E(\widebar{\rho}_G)|}}{\lambda^{1 + (r-1)n }c_\lambda^{n|E(G)|}}\cdot 
  \frac{\lambda^{|V(\widehat{\rho}_G)|}c_\lambda^{|E(\widehat{\rho}_G)|}}{\lambda^{r}c_\lambda^{|E(G)|}}\cdot
  \frac{\lambda^{-|V(H)|}c_\lambda^{-|E(H)|}}{\lambda^{-1}}
  \\
  \nonumber
  &=&\frac{\lambda^{|V(\widebar{\rho}_G)|}c_\lambda^{|E(\widebar{\rho}_G)|}}{\lambda^{1 + (r-1)n }c_\lambda^{n|E(G)|}}\cdot\big(
  \lambda^{|V(H)|-1} c_{\lambda}^{|E(H)|}\big)^{-1}  \\
  \nonumber
  &=& O(1), 
\end{eqnarray}
 therefore \eqref{negli} holds at the rank $n+1$. 
 As a consequence, the leading terms in \eqref{cumulant-diagram1} are those associated
  with the connected partition diagrams $\Gamma(\rho ,\pi )$ having the highest
  block count, i.e. which have $1+(r-1)n$ blocks, see Figure~\ref{fig:diagram3}
  for a sample of such a partition diagram.

\begin{figure}[H]
\captionsetup[subfigure]{font=footnotesize}
\centering
\subcaptionbox{Diagram $\Gamma(\rho,\pi)$ and graph $\widetilde{\rho}_G$ in blue.}[.49\textwidth]{
\begin{tikzpicture}[scale=0.9] 
\draw[black, thick] (0,0) rectangle (5,6);

\node[anchor=east,font=\small] at (0.8,5) {1};
\node[anchor=east,font=\small] at (0.8,4) {2};
\node[anchor=east,font=\small] at (0.8,3) {3};
\node[anchor=east,font=\small] at (0.8,2) {4};
\node[anchor=east,font=\small] at (0.8,1) {5};

\node[anchor=south,font=\small] at (1,0) {1};
\node[anchor=south,font=\small] at (2,0) {2};
\node[anchor=south,font=\small] at (3,0) {3};
\node[anchor=south,font=\small] at (4,0) {4};

\filldraw [gray] (1,1) circle (2pt);
\filldraw [gray] (2,1) circle (2pt);
\filldraw [gray] (3,1) circle (2pt);
\filldraw [gray] (4,1) circle (2pt);
\filldraw [gray] (1,2) circle (2pt);
\filldraw [gray] (2,2) circle (2pt);
\filldraw [gray] (3,2) circle (2pt);
\filldraw [gray] (4,2) circle (2pt);
\filldraw [gray] (1,3) circle (2pt);
\filldraw [gray] (2,3) circle (2pt);
\filldraw [gray] (3,3) circle (2pt);
\filldraw [gray] (4,3) circle (2pt);
\filldraw [gray] (2,3) circle (2pt);
\filldraw [gray] (1,4) circle (2pt);
\filldraw [gray] (2,4) circle (2pt);
\filldraw [gray] (3,4) circle (2pt);
\filldraw [gray] (4,4) circle (2pt);
\filldraw [gray] (1,5) circle (2pt);
\filldraw [gray] (2,5) circle (2pt);
\filldraw [gray] (3,5) circle (2pt);
\filldraw [gray] (4,5) circle (2pt);

\draw[very thick] (1,5) -- (1,1); 

\foreach \i in {1,...,5}
         {
           \draw[thick,dash dot,blue] (1,\i) .. controls (1.5,\i+.5) .. (2,\i);
           \draw[thick,dash dot,blue] (2,\i) .. controls (2.5,\i+.5) .. (3,\i);
           \draw[thick,dash dot,blue] (2,\i) .. controls (3,\i-.5) .. (4,\i);
           \draw[thick,dash dot,blue] (3,\i) .. controls (3.5,\i+.5) .. (4,\i);
         }
         
  \begin{pgfonlayer}{background}
    \filldraw [line width=4mm,black!3]
      (0.2,0.2)  rectangle (4.8,5.8);
  \end{pgfonlayer}
\end{tikzpicture}}
\subcaptionbox{Diagram $\Gamma(\rho ,\pi)$ and graph $\rho_G$ in red.}[.49\textwidth]{
\begin{tikzpicture}[scale=0.9] 
\draw[black, thick] (0,0) rectangle (5,6);

\node[anchor=east,font=\small] at (0.8,5) {1};
\node[anchor=east,font=\small] at (0.8,4) {2};
\node[anchor=east,font=\small] at (0.8,3) {3};
\node[anchor=east,font=\small] at (0.8,2) {4};
\node[anchor=east,font=\small] at (0.8,1) {5};

\node[anchor=south,font=\small] at (1,0) {1};
\node[anchor=south,font=\small] at (2,0) {2};
\node[anchor=south,font=\small] at (3,0) {3};
\node[anchor=south,font=\small] at (4,0) {4};

\filldraw [gray] (1,1) circle (2pt);
\filldraw [gray] (2,1) circle (2pt);
\filldraw [gray] (3,1) circle (2pt);
\filldraw [gray] (4,1) circle (2pt);
\filldraw [gray] (1,2) circle (2pt);
\filldraw [gray] (2,2) circle (2pt);
\filldraw [gray] (3,2) circle (2pt);
\filldraw [gray] (4,2) circle (2pt);
\filldraw [gray] (1,3) circle (2pt);
\filldraw [gray] (2,3) circle (2pt);
\filldraw [gray] (3,3) circle (2pt);
\filldraw [gray] (4,3) circle (2pt);
\filldraw [gray] (2,3) circle (2pt);
\filldraw [gray] (1,4) circle (2pt);
\filldraw [gray] (2,4) circle (2pt);
\filldraw [gray] (3,4) circle (2pt);
\filldraw [gray] (4,4) circle (2pt);
\filldraw [gray] (1,5) circle (2pt);
\filldraw [gray] (2,5) circle (2pt);
\filldraw [gray] (3,5) circle (2pt);
\filldraw [gray] (4,5) circle (2pt);

\draw[very thick] (1,5) -- (1,1); 

\foreach \i in {1,...,5}
         {
           \draw[thick,dash dot,purple] (1,\i) .. controls (1.5,\i+.5) .. (2,\i);
           \draw[thick,dash dot,purple] (2,\i) .. controls (2.5,\i+.5) .. (3,\i);
           \draw[thick,dash dot,purple] (2,\i) .. controls (3,\i-.5) .. (4,\i);
           \draw[thick,dash dot,purple] (3,\i) .. controls (3.5,\i+.5) .. (4,\i);
         }
  \begin{pgfonlayer}{background}
    \filldraw [line width=4mm,black!3]
      (0.2,0.2)  rectangle (4.8,5.8);
  \end{pgfonlayer}
\end{tikzpicture}}
\caption{Example of maximal connected partition diagram with $n=5$ and $r=4$.}
\label{fig:diagram3}
\end{figure}

\vspace{-0.5cm} 

\noindent
 Finally, we observe that any maximal connected non-flat partition
 $\rho \in {\rm M}(n,r)$ satisfies
 $|E(\rho_G)|=n\times|E(G)|$, 
as can be checked in Figure~\ref{fig:diagram3}.
 Therefore, by 
 \eqref{coeff-0}-\eqref{coeff-1},
 \eqref{cumulant-diagram1} and \eqref{integ-connecting3}, we obtain 
\begin{eqnarray*} 
  \lefteqn{
    c^{n|E(G)|}
    C^{1+(r-1)n}
    c_\lambda^{n|E(G)|}
    ( (r-1)r )^{n-1}(n-1)!
    \lambda^{1+(r-1)n}
  }
  \\
   & \leq &   
    \lambda^{1+(r-1)n}
  \sum_{\rho\in {\rm M}(n,r)}\int_{(\R^d)^{1+(r-1)n}}
  \Bigg(
  \prod_{\{\eta_1,\eta_2\}\in E(\rho_G)}H_\lambda(x_{\eta_1},x_{\eta_2})
  \Bigg)
  \prod_{\eta\in V(\rho_G)}\mu(\mathrm{d}x_{\eta}),
  \\
   & \leq &     \kappa_n(N_G)
      \\
  & \leq & 
  n!^r r!^{n-1} 
  ( 1 + \mu ( \real^d))^{1+(r-1)n}
  c_\lambda^{n|E(G)|}
 \lambda^{1+(r-1)n}, 
\end{eqnarray*}
which yields \eqref{equiv-1}.
\end{Proof}
 In what follows, we consider the centered and normalized subgraph count cumulants defined as 
 $$
 \widetilde{N}_G
 := \frac{N_G - \kappa_1 (N_G)}{\sqrt{\kappa_2(N_G)}}, \qquad n \geq 1. 
$$
 The following result shows that for $n\geq 3$ the normalized cumulant
 $\kappa_n(\widetilde{N}_G)$ tends to zero in \eqref{Statuleviciuscond},
 hence $\widetilde{N}_G$ converges in distribution to the normal
 distribution by Theorem~1 in \cite{Janson1988}.   
\begin{corollary}[Dilute regime]
\label{t1-c}
  Let $G$ be a connected graph with $|V(G)| = r$ vertices, $r\geq 2$, 
satisfying Assumption~\ref{a61} for $n=2$ in the dilute regime \eqref{fjnldsf}. 
  We have the normalized cumulant bounds 
 \begin{equation}
    \label{Statuleviciuscond}
  \big|\kappa_n \big(\widetilde{N}_G\big)\big|
  \leq n!^r ( K \lambda )^{-(n/2-1)},
  \qquad \lambda \geq 1, \quad n\geq 2,
\end{equation}
where $K>0$ is a constant independent of $\lambda >0$ and $n\geq 1$.
\end{corollary}
\begin{Proof}
  We note that the upper bound in \eqref{equiv-1}
  does not require Assumption~\ref{a61},
  hence we have, for $n\geq 2$,
$$ 
   \big|\kappa_n(\widetilde{N}_G)\big|
  \leq 
   \frac{n!^r 
     c_\lambda^{n|E(G)|} (K_2 \lambda )^{1+(r-1)n}}{\big((2-1)! c_\lambda^{2|E(G)|}
   (K_1 \lambda )^{1+2(r-1)}\big)^{n/2}}
 = 
   K_2
   \left(
   \frac{ 
      (K_2/K_1)^{r-1}}{\sqrt{K_1}}
   \right)^n
   n!^r \lambda^{-(n/2-1)}. 
$$ 
\end{Proof}
The following result yields a positive cumulant growth
of order $\alpha     -(\alpha - 1)r>0$
in \eqref{cumulant-rhop2-0} for trees in the sparse regime
with $\alpha \in [1, r/(r-1) )$,
while in the case of non-tree graphs such as
cycle graphs the growth rate {exponent} 
$r - \alpha |E(G)|
\leq ( 1 - \alpha ) r \leq 0$
is negative or zero in \eqref{cumulant-rhoph1} and \eqref{cumulant-rhoph2}. 
In addition, the normalized cumulant 
$\kappa_n(\widetilde{N}_G)$ tends to zero for $n\geq 3$ in \eqref{jfkla} only
when $G$ is a tree, in which case 
$\widetilde{N}_G$ converges in distribution to the normal
distribution by Theorem~1 in \cite{Janson1988}.   
We note that when $\alpha = 1$, \eqref{jfkla} is consistent with
 \eqref{Statuleviciuscond}. 
\begin{prop}[Sparse regime] 
\label{th6.4}
    Let $G$ be a connected graph with $|V(G)|=r$ vertices, $r\geq 2$,
  satisfying Assumption~\ref{a61} for all $n\geq 1$
  in the sparse regime \eqref{fjnldsf-2}. 
  \begin{enumerate}[a)] 
  \item 
 If $G$ is a tree, i.e. $|E(G)| = r-1$, we have the cumulant bounds 
\begin{equation} 
\label{cumulant-rhop2-0}
  (K_1)^r
 \lambda^{\alpha  -(\alpha - 1)r }
     \leq 
  \kappa_n(N_G)
    \leq 
  n!^r
    (K_2)^r \lambda^{ \alpha     -(\alpha - 1)r       } ,
  \quad \lambda \geq 1, 
\end{equation} 
 for some constants $K_1>0$, $K_2>1$ independent of $\lambda, n\geq 1$.
\item
 If $G$ is not a tree, i.e. $|E(G)|\geq r$,
 we have the cumulant bounds 
 \begin{equation} 
 \label{cumulant-rhoph1}
  (K_1)^r 
 \lambda^{r-\alpha |E(G)|}
 \leq 
  \kappa_n(N_G)
  \leq
   n!^r
     (K_2)^r
   \lambda^{r-\alpha |E(G)|}, 
   \quad \lambda \geq 1,
 \end{equation} 
 for some constants $K_1>0$, $K_2>1$ independent of $\lambda, n\geq 1$. 
\item
  If $G$ is a cycle, i.e. $|E(G)| = r$,
  we have the cumulant bounds 
 \begin{equation} 
\label{cumulant-rhoph2}
  (K_1)^r 
 \lambda^{- (\alpha - 1 )r}
\leq 
  \kappa_n(N_G)
  \leq
  n!^r
   (K_2)^r
  \lambda^{ - (\alpha -1)r},
    \quad \lambda \geq 1, 
\end{equation} 
 for some constants $K_1>0$, $K_2>1$ independent of $\lambda, n\geq 1$. 
\end{enumerate}
\end{prop}
\begin{Proof}
  In the sparse regime \eqref{fjnldsf-2}, 
  every edge in the graph $\rho_G$ contributes a power $\lambda^{-\alpha}$
  and every vertex contributes a power $\lambda$,
  hence every term in \eqref{cumulant-diagram1} contributes a power
  \begin{equation}
    \label{fjkld34} 
  \lambda^{|V(\rho_G)| - \alpha |E(\rho_G)|}
  =
  \lambda^{    \alpha    - (\alpha - 1 ) |V(\rho_G)|    + ( |V (\rho_G)| - |E(\rho_G)| - 1 ) \alpha  }
 \leq
 \lambda^{   \alpha   - (\alpha - 1 ) |V(\rho_G)|     }
\end{equation} 
  since $|V (\rho_G)| - |E(\rho_G)| - 1 \leq 0$. 
  In addition, for any connected partition 
 $\rho\in\Pi_{\widehat{1}} ([n]\times[r])$,
  we have 
  $$r\le|V(\rho_G)|\leq 1+(r-1)n. 
  $$ 
    
\noindent
 $a)$
 When $G$ is a tree and the graph $\rho_G$ is also a tree, 
 i.e. $|V (\rho_G)| - |E(\rho_G)| - 1 = 0$, and 
 the maximal order 
 $\lambda^{
   \alpha
   - (\alpha - 1 ) |V(\rho_G)|
     }
 $ is attained in \eqref{fjkld34}, 
  see Figure~\ref{fig:diagram2-1-4} for an example.

\smallskip

\begin{figure}[H]
\captionsetup[subfigure]{font=footnotesize}
\centering
\subcaptionbox{Diagram $\Gamma(\rho,\pi)$ and multigraph $\widetilde{\rho}_G$ in blue.}[.49\textwidth]{
\begin{tikzpicture}[scale=0.9] 
\draw[black, thick] (0,0) rectangle (5,6);

\node[anchor=east,font=\small] at (0.8,5) {1};
\node[anchor=east,font=\small] at (0.8,4) {2};
\node[anchor=east,font=\small] at (0.8,3) {3};
\node[anchor=east,font=\small] at (0.8,2) {4};
\node[anchor=east,font=\small] at (0.8,1) {5};

\node[anchor=south,font=\small] at (1,0) {1};
\node[anchor=south,font=\small] at (2,0) {2};
\node[anchor=south,font=\small] at (3,0) {3};
\node[anchor=south,font=\small] at (4,0) {4};

\filldraw [gray] (1,1) circle (2pt);
\filldraw [gray] (2,1) circle (2pt);
\filldraw [gray] (3,1) circle (2pt);
\filldraw [gray] (4,1) circle (2pt);
\filldraw [gray] (1,2) circle (2pt);
\filldraw [gray] (2,2) circle (2pt);
\filldraw [gray] (3,2) circle (2pt);
\filldraw [gray] (4,2) circle (2pt);
\filldraw [gray] (1,3) circle (2pt);
\filldraw [gray] (2,3) circle (2pt);
\filldraw [gray] (3,3) circle (2pt);
\filldraw [gray] (4,3) circle (2pt);
\filldraw [gray] (2,3) circle (2pt);
\filldraw [gray] (1,4) circle (2pt);
\filldraw [gray] (2,4) circle (2pt);
\filldraw [gray] (3,4) circle (2pt);
\filldraw [gray] (4,4) circle (2pt);
\filldraw [gray] (1,5) circle (2pt);
\filldraw [gray] (2,5) circle (2pt);
\filldraw [gray] (3,5) circle (2pt);
\filldraw [gray] (4,5) circle (2pt);

\draw[very thick] (1,5) -- (1,4); 
\draw[very thick] (3,5) -- (4,4);

\draw[very thick] (1,2) -- (1,1);
\draw[very thick] (2,2) -- (2,4);
\draw[very thick] (2,1) -- (3,2) -- (4,3) -- (3,4);

\foreach \i in {1,...,5}
         {
           \draw[thick,dash dot,blue] (1,\i) .. controls (1.5,\i+.5) .. (2,\i);
           \draw[thick,dash dot,blue] (2,\i) .. controls (2.5,\i+.5) .. (3,\i);
           \draw[thick,dash dot,blue] (2,\i) .. controls (3,\i-.5) .. (4,\i);
         }
         
  \begin{pgfonlayer}{background}
    \filldraw [line width=4mm,black!3]
      (0.2,0.2)  rectangle (4.8,5.8);
  \end{pgfonlayer}
\end{tikzpicture}}
\subcaptionbox{Diagram $\Gamma(\sigma,\pi)$ and graph $\rho_G$ in red.}[.49\textwidth]{
\begin{tikzpicture}[scale=0.9] 
\draw[black, thick] (0,0) rectangle (5,6);

\node[anchor=east,font=\small] at (0.8,5) {1};
\node[anchor=east,font=\small] at (0.8,4) {2};
\node[anchor=east,font=\small] at (0.8,3) {3};
\node[anchor=east,font=\small] at (0.8,2) {4};
\node[anchor=east,font=\small] at (0.8,1) {5};

\node[anchor=south,font=\small] at (1,0) {1};
\node[anchor=south,font=\small] at (2,0) {2};
\node[anchor=south,font=\small] at (3,0) {3};
\node[anchor=south,font=\small] at (4,0) {4};

\filldraw [gray] (1,1) circle (2pt);
\filldraw [gray] (2,1) circle (2pt);
\filldraw [gray] (3,1) circle (2pt);
\filldraw [gray] (4,1) circle (2pt);
\filldraw [gray] (1,2) circle (2pt);
\filldraw [gray] (2,2) circle (2pt);
\filldraw [gray] (3,2) circle (2pt);
\filldraw [gray] (4,2) circle (2pt);
\filldraw [gray] (1,3) circle (2pt);
\filldraw [gray] (2,3) circle (2pt);
\filldraw [gray] (3,3) circle (2pt);
\filldraw [gray] (4,3) circle (2pt);
\filldraw [gray] (2,3) circle (2pt);
\filldraw [gray] (1,4) circle (2pt);
\filldraw [gray] (2,4) circle (2pt);
\filldraw [gray] (3,4) circle (2pt);
\filldraw [gray] (4,4) circle (2pt);
\filldraw [gray] (1,5) circle (2pt);
\filldraw [gray] (2,5) circle (2pt);
\filldraw [gray] (3,5) circle (2pt);
\filldraw [gray] (4,5) circle (2pt);

\draw[very thick] (1,5) -- (1,4); 
\draw[very thick] (3,5) -- (4,4);

\draw[very thick] (1,2) -- (1,1);
\draw[very thick] (2,2) -- (2,4);
\draw[very thick] (2,1) -- (3,2) -- (4,3) -- (3,4);

\foreach \i in {1,...,3}
         {
           \draw[thick,dash dot,purple] (1,\i) .. controls (1.5,\i+.5) .. (2,\i);
           \draw[thick,dash dot,purple] (2,\i) .. controls (2.5,\i+.5) .. (3,\i);
           \draw[thick,dash dot,purple] (2,\i) .. controls (3,\i-.5) .. (4,\i);
         }

           \draw[thick,dash dot,purple] (1,4) .. controls (1.5,4+.5) .. (2,4);
           \draw[thick,dash dot,purple] (2,4) .. controls (3,4-.5) .. (4,4);
\foreach \i in {5,...,5}
         {
           \draw[thick,dash dot,purple] (1,\i) .. controls (1.5,\i+.5) .. (2,\i);
           \draw[thick,dash dot,purple] (2,\i) .. controls (2.5,\i+.5) .. (3,\i);
                      \draw[thick,dash dot,purple] (2,\i) .. controls (3,\i-.5) .. (4,\i);
         }
  \begin{pgfonlayer}{background}
    \filldraw [line width=4mm,black!3]
      (0.2,0.2)  rectangle (4.8,5.8);
  \end{pgfonlayer}
\end{tikzpicture}}
\caption{Example of connected partition diagram with $\rho_G$ a tree and $n=5$, $r=4$.}
\label{fig:diagram2-1-4}
\end{figure}

\vskip-0.3cm

\noindent
 In this case, the corresponding term in \eqref{cumulant-diagram1} contributes a power
$$
 \lambda^{|V(\rho_G)|-\alpha|E(\rho_G)|}
  =
  \lambda^{
    \alpha
    - (\alpha - 1 ) |V(\rho_G)|
  }, 
  \qquad \lambda \geq 1. 
  $$
 In this case, since $|V(\rho_G)|\geq r$ and $\alpha \geq 1$, the optimal rate
 $\lambda^{ \alpha - ( \alpha - 1 ) r }$
 is attained by the partition diagrams
  $\Gamma(\rho ,\pi )$ such that 
  $|V(\rho_G )|=r$, as illustrated in Figure~\ref{fig:diagram4-0}. 
  
\begin{figure}[H]
\captionsetup[subfigure]{font=footnotesize}
\centering
\subcaptionbox{Diagram $\Gamma(\rho,\pi)$ and multigraph $\widetilde{\rho}_G$ in blue.}[.49\textwidth]{
\begin{tikzpicture}[scale=0.9] 

\draw[black, thick] (0,0) rectangle (5,6);

\node[anchor=east,font=\small] at (0.8,5) {1};
\node[anchor=east,font=\small] at (0.8,4) {2};
\node[anchor=east,font=\small] at (0.8,3) {3};
\node[anchor=east,font=\small] at (0.8,2) {4};
\node[anchor=east,font=\small] at (0.8,1) {5};

\node[anchor=south,font=\small] at (1,0) {1};
\node[anchor=south,font=\small] at (2,0) {2};
\node[anchor=south,font=\small] at (3,0) {3};
\node[anchor=south,font=\small] at (4,0) {4};

\filldraw [gray] (1,1) circle (2pt);
\filldraw [gray] (2,1) circle (2pt);
\filldraw [gray] (3,1) circle (2pt);
\filldraw [gray] (4,1) circle (2pt);
\filldraw [gray] (1,2) circle (2pt);
\filldraw [gray] (2,2) circle (2pt);
\filldraw [gray] (3,2) circle (2pt);
\filldraw [gray] (4,2) circle (2pt);
\filldraw [gray] (1,3) circle (2pt);
\filldraw [gray] (2,3) circle (2pt);
\filldraw [gray] (3,3) circle (2pt);
\filldraw [gray] (4,3) circle (2pt);
\filldraw [gray] (2,3) circle (2pt);
\filldraw [gray] (1,4) circle (2pt);
\filldraw [gray] (2,4) circle (2pt);
\filldraw [gray] (3,4) circle (2pt);
\filldraw [gray] (4,4) circle (2pt);
\filldraw [gray] (1,5) circle (2pt);
\filldraw [gray] (2,5) circle (2pt);
\filldraw [gray] (3,5) circle (2pt);
\filldraw [gray] (4,5) circle (2pt);

\foreach \i in {1,...,4}
{
\draw[very thick] (\i,5) -- (\i,1); 
}

  \foreach \i in {1,...,5}
         {
           \draw[thick,dash dot,blue] (1,\i) .. controls (1.5,\i+.5) .. (2,\i);
           \draw[thick,dash dot,blue] (2,\i) .. controls (2.5,\i+.5) .. (3,\i);
           \draw[thick,dash dot,blue] (2,\i) .. controls (3,\i-.5) .. (4,\i);
         }
         
  \begin{pgfonlayer}{background}
    \filldraw [line width=4mm,black!3]
      (0.2,0.2)  rectangle (4.8,5.8);
  \end{pgfonlayer}
\end{tikzpicture}}
\subcaptionbox{Diagram $\Gamma(\rho,\pi)$ and graph $\rho_G$ in red.}[.49\textwidth]{
\begin{tikzpicture}[scale=0.9] 
\draw[black, thick] (0,0) rectangle (5,6);

\node[anchor=east,font=\small] at (0.8,5) {1};
\node[anchor=east,font=\small] at (0.8,4) {2};
\node[anchor=east,font=\small] at (0.8,3) {3};
\node[anchor=east,font=\small] at (0.8,2) {4};
\node[anchor=east,font=\small] at (0.8,1) {5};

\node[anchor=south,font=\small] at (1,0) {1};
\node[anchor=south,font=\small] at (2,0) {2};
\node[anchor=south,font=\small] at (3,0) {3};
\node[anchor=south,font=\small] at (4,0) {4};

\filldraw [gray] (1,1) circle (2pt);
\filldraw [gray] (2,1) circle (2pt);
\filldraw [gray] (3,1) circle (2pt);
\filldraw [gray] (4,1) circle (2pt);
\filldraw [gray] (1,2) circle (2pt);
\filldraw [gray] (2,2) circle (2pt);
\filldraw [gray] (3,2) circle (2pt);
\filldraw [gray] (4,2) circle (2pt);
\filldraw [gray] (1,3) circle (2pt);
\filldraw [gray] (2,3) circle (2pt);
\filldraw [gray] (3,3) circle (2pt);
\filldraw [gray] (4,3) circle (2pt);
\filldraw [gray] (2,3) circle (2pt);
\filldraw [gray] (1,4) circle (2pt);
\filldraw [gray] (2,4) circle (2pt);
\filldraw [gray] (3,4) circle (2pt);
\filldraw [gray] (4,4) circle (2pt);
\filldraw [gray] (1,5) circle (2pt);
\filldraw [gray] (2,5) circle (2pt);
\filldraw [gray] (3,5) circle (2pt);
\filldraw [gray] (4,5) circle (2pt);

\foreach \i in {1,...,4}
{
\draw[very thick] (\i,5) -- (\i,1); 
}

\foreach \i in {5,...,5}
         {
           \draw[thick,dash dot,purple] (1,\i) .. controls (1.5,\i+.5) .. (2,\i);
           \draw[thick,dash dot,purple] (2,\i) .. controls (2.5,\i+.5) .. (3,\i);
                      \draw[thick,dash dot,purple] (2,\i) .. controls (3,\i-.5) .. (4,\i);
         }
         
  \begin{pgfonlayer}{background}
    \filldraw [line width=4mm,black!3]
      (0.2,0.2)  rectangle (4.8,5.8);
  \end{pgfonlayer}
\end{tikzpicture}}
\caption{Tree diagram $\rho_G$ with $G$ a tree with $|V(\rho_G)|=r$ and $n=5$, $r=4$.}
\label{fig:diagram4-0}
\end{figure}

\vskip-0.3cm

\noindent
 We conclude to \eqref{cumulant-rhop2-0} 
 as in the proof of Proposition~\ref{t1}, 
 by upper bounding the count of connected non-flat
 partitions from \eqref{coeff-0}
 and by lower bounding it by $1$. 

 \noindent
   $b)$ 
   When $G$ is not a tree it contains at least one cycle, and
   for any partition $\rho\in\Pi_{\widehat{1}} ([n]\times[r])$  
   the same holds for the graph $\rho_G$. 
In this case, the highest order contribution in
\eqref{cumulant-diagram1} is attained by connected non-flat partition diagrams 
$\Gamma ( \rho , \pi )$, $\rho\in\Pi_{\widehat{1}} ([n]\times[r])$,
such that $\rho_G$ has $|V(\rho_G)| = r$ vertices, 
and their contribution is given by a power of order $\lambda^{r-\alpha |E(G)|}$. 
 An example of such partition $\rho$ is given in Figure~\ref{fig:diagram4}, 
 with $G$ a cycle. 

\smallskip

\begin{figure}[H]
\captionsetup[subfigure]{font=footnotesize}
\centering
\subcaptionbox{Diagram $\Gamma(\rho,\pi)$ and multigraph $\widetilde{\rho}_G$ in blue.}[.49\textwidth]{
\begin{tikzpicture}[scale=0.9] 
\draw[black, thick] (0,0) rectangle (5,4);

\node[anchor=east,font=\small] at (0.8,3) {1};
\node[anchor=east,font=\small] at (0.8,2) {2};
\node[anchor=east,font=\small] at (0.8,1) {3};

\node[anchor=south,font=\small] at (1,0) {1};
\node[anchor=south,font=\small] at (2,0) {2};
\node[anchor=south,font=\small] at (3,0) {3};
\node[anchor=south,font=\small] at (4,0) {4};

\filldraw [gray] (1,1) circle (2pt);
\filldraw [gray] (2,1) circle (2pt);
\filldraw [gray] (3,1) circle (2pt);
\filldraw [gray] (4,1) circle (2pt);
\filldraw [gray] (1,2) circle (2pt);
\filldraw [gray] (2,2) circle (2pt);
\filldraw [gray] (3,2) circle (2pt);
\filldraw [gray] (4,2) circle (2pt);
\filldraw [gray] (1,3) circle (2pt);
\filldraw [gray] (2,3) circle (2pt);
\filldraw [gray] (3,3) circle (2pt);
\filldraw [gray] (4,3) circle (2pt);
\filldraw [gray] (2,3) circle (2pt);

\foreach \i in {1,...,4}
{
\draw[very thick] (\i,3) -- (\i,1); 
}

\foreach \i in {1,...,3}
         {
           \draw[thick,dash dot,blue] (1,\i) .. controls (1.5,\i+.5) .. (2,\i);
           \draw[thick,dash dot,blue] (2,\i) .. controls (2.5,\i+.5) .. (3,\i);
           \draw[thick,dash dot,blue] (3,\i) .. controls (3.5,\i+.5) .. (4,\i);
           \draw[thick,dash dot,blue] (4,\i) .. controls (2.5,\i-.5) .. (1,\i);
         }
         
  \begin{pgfonlayer}{background}
    \filldraw [line width=4mm,black!3]
      (0.2,0.2)  rectangle (4.8,3.8);
  \end{pgfonlayer}
\end{tikzpicture}}
\subcaptionbox{Diagram $\Gamma(\sigma,\pi)$ and graph $\rho_G$ in red.}[.49\textwidth]{
\begin{tikzpicture}[scale=0.9] 
\draw[black, thick] (0,0) rectangle (5,4);

\node[anchor=east,font=\small] at (0.8,3) {1};
\node[anchor=east,font=\small] at (0.8,2) {2};
\node[anchor=east,font=\small] at (0.8,1) {3};

\node[anchor=south,font=\small] at (1,0) {1};
\node[anchor=south,font=\small] at (2,0) {2};
\node[anchor=south,font=\small] at (3,0) {3};
\node[anchor=south,font=\small] at (4,0) {4};

\filldraw [gray] (1,1) circle (2pt);
\filldraw [gray] (2,1) circle (2pt);
\filldraw [gray] (3,1) circle (2pt);
\filldraw [gray] (4,1) circle (2pt);
\filldraw [gray] (1,2) circle (2pt);
\filldraw [gray] (2,2) circle (2pt);
\filldraw [gray] (3,2) circle (2pt);
\filldraw [gray] (4,2) circle (2pt);
\filldraw [gray] (1,3) circle (2pt);
\filldraw [gray] (2,3) circle (2pt);
\filldraw [gray] (3,3) circle (2pt);
\filldraw [gray] (4,3) circle (2pt);
\filldraw [gray] (2,3) circle (2pt);

\foreach \i in {1,...,4}
{
\draw[very thick] (\i,3) -- (\i,1); 
}

           \draw[thick,dash dot,purple] (1,3) .. controls (1.5,3+.5) .. (2,3);
           \draw[thick,dash dot,purple] (2,3) .. controls (2.5,3+.5) .. (3,3);
           \draw[thick,dash dot,purple] (3,3) .. controls (3.5,3+.5) .. (4,3);
           \draw[thick,dash dot,purple] (4,3) .. controls (2.5,3-.5) .. (1,3);
         
  \begin{pgfonlayer}{background}
    \filldraw [line width=4mm,black!3]
      (0.2,0.2)  rectangle (4.8,3.8);
  \end{pgfonlayer}
\end{tikzpicture}}
\caption{Cycle graph $\rho_G$ with $G$ a cycle graph and $n=5$, $r=4$.}
\label{fig:diagram4}
\end{figure}

\vskip-0.3cm

\noindent
 Indeed, in order to remain non-flat,
 the partition $\rho$ can only be modified
 into a partition $\sigma$ by splitting a block of $\rho_G$ in two, 
 which entails the addition of a number $q$ of edges, $q\geq 1$, 
 resulting into an additional factor $\lambda^{1- q \alpha} \leq 1$ 
 that may only lower the order of the contribution, see
 Figures~\ref{fig:diagram3-1-0}-\ref{fig:diagram3-1-3}
 for examples with $G$ a graph containing one cycle.
 
\begin{figure}[H]
\captionsetup[subfigure]{font=footnotesize}
\centering
\subcaptionbox{Diagram $\Gamma(\rho,\pi)$ with order $\lambda^{4-4\alpha}$.}[.47\textwidth]{\begin{tikzpicture}[scale=0.9] 
\draw[black, thick] (0,0) rectangle (5,4);

\node[anchor=east,font=\small] at (0.8,3) {1};
\node[anchor=east,font=\small] at (0.8,2) {2};
\node[anchor=east,font=\small] at (0.8,1) {3};
\node[anchor=south,font=\small] at (1,0) {1};
\node[anchor=south,font=\small] at (2,0) {2};
\node[anchor=south,font=\small] at (3,0) {3};
\node[anchor=south,font=\small] at (4,0) {4};

\filldraw [gray] (1,1) circle (2pt);
\filldraw [gray] (2,1) circle (2pt);
\filldraw [gray] (3,1) circle (2pt);
\filldraw [gray] (4,1) circle (2pt);
\filldraw [gray] (1,2) circle (2pt);
\filldraw [gray] (2,2) circle (2pt);
\filldraw [gray] (3,2) circle (2pt);
\filldraw [gray] (4,2) circle (2pt);
\filldraw [gray] (1,3) circle (2pt);
\filldraw [gray] (2,3) circle (2pt);
\filldraw [gray] (3,3) circle (2pt);
\filldraw [gray] (4,3) circle (2pt);
\filldraw [gray] (2,3) circle (2pt);

\draw[very thick] (1,1) -- (1,1) -- (1,3);
\draw[very thick] (2,1) -- (2,1) -- (2,3);
\draw[very thick] (3,1) -- (3,1) -- (3,3);
\draw[very thick] (4,1) -- (4,1) -- (4,3);
\foreach \i in {3}
         {
           \draw[thick,dash dot,purple] (1,\i) .. controls (1.5,\i+.5) .. (2,\i);
           \draw[thick,dash dot,purple] (2,\i) .. controls (2.5,\i+.5) .. (3,\i);
           \draw[thick,dash dot,purple] (3,\i) .. controls (3.5,\i+.5) .. (4,\i);
           \draw[thick,dash dot,purple] (2,\i) .. controls (3,\i-.5) .. (4,\i);
         }
         
  \begin{pgfonlayer}{background}
    \filldraw [line width=4mm,black!3]
      (0.2,0.2)  rectangle (4.8,3.8);
  \end{pgfonlayer}
\end{tikzpicture}}
\hfill
\subcaptionbox{Diagram $\Gamma(\sigma,\pi)$ with order $\lambda^{5-5\alpha}= \lambda^{4-4\alpha}\lambda^{- ( \alpha - 1)}$.}[.52\textwidth]{
\begin{tikzpicture}[scale=0.9] 
\draw[black, thick] (0,0) rectangle (5,4);

\node[anchor=east,font=\small] at (0.8,3) {1};
\node[anchor=east,font=\small] at (0.8,2) {2};
\node[anchor=east,font=\small] at (0.8,1) {3};
\node[anchor=south,font=\small] at (1,0) {1};
\node[anchor=south,font=\small] at (2,0) {2};
\node[anchor=south,font=\small] at (3,0) {3};
\node[anchor=south,font=\small] at (4,0) {4};

\filldraw [gray] (1,1) circle (2pt);
\filldraw [gray] (2,1) circle (2pt);
\filldraw [gray] (3,1) circle (2pt);
\filldraw [gray] (4,1) circle (2pt);
\filldraw [gray] (1,2) circle (2pt);
\filldraw [gray] (2,2) circle (2pt);
\filldraw [gray] (3,2) circle (2pt);
\filldraw [gray] (4,2) circle (2pt);
\filldraw [gray] (1,3) circle (2pt);
\filldraw [gray] (2,3) circle (2pt);
\filldraw [gray] (3,3) circle (2pt);
\filldraw [gray] (4,3) circle (2pt);
\filldraw [gray] (2,3) circle (2pt);

\draw[very thick] (1,1) -- (1,2);
\draw[very thick] (2,1) -- (2,1) -- (2,3);
\draw[very thick] (3,1) -- (3,1) -- (3,3);
\draw[very thick] (4,1) -- (4,1) -- (4,3);

\draw[thick,dash dot,purple] (1,2) .. controls (1.5,2.5) .. (2,2);

\foreach \i in {3}
         {
           \draw[thick,dash dot,purple] (1,\i) .. controls (1.5,\i+.5) .. (2,\i);
           \draw[thick,dash dot,purple] (2,\i) .. controls (2.5,\i+.5) .. (3,\i);
           \draw[thick,dash dot,purple] (3,\i) .. controls (3.5,\i+.5) .. (4,\i);
           \draw[thick,dash dot,purple] (2,\i) .. controls (3,\i-.5) .. (4,\i);
         }
         
  \begin{pgfonlayer}{background}
    \filldraw [line width=4mm,black!3]
      (0.2,0.2)  rectangle (4.8,3.8);
  \end{pgfonlayer}
\end{tikzpicture}}
\caption{Splitting of a vertex with addition of one edge and $n=3$, $r=4$.}
\label{fig:diagram3-1-0}
\end{figure}

\vskip-0.3cm

\noindent

\smallskip

\begin{figure}[H]
\captionsetup[subfigure]{font=footnotesize}
\centering
\subcaptionbox{Diagram $\Gamma(\rho,\pi)$ with order $\lambda^{4-4\alpha}$.}[.49\textwidth]{
\begin{tikzpicture}[scale=0.9] 
\draw[black, thick] (0,0) rectangle (5,4);

\node[anchor=east,font=\small] at (0.8,3) {1};
\node[anchor=east,font=\small] at (0.8,2) {2};
\node[anchor=east,font=\small] at (0.8,1) {3};
\node[anchor=south,font=\small] at (1,0) {1};
\node[anchor=south,font=\small] at (2,0) {2};
\node[anchor=south,font=\small] at (3,0) {3};
\node[anchor=south,font=\small] at (4,0) {4};

\filldraw [gray] (1,1) circle (2pt);
\filldraw [gray] (2,1) circle (2pt);
\filldraw [gray] (3,1) circle (2pt);
\filldraw [gray] (4,1) circle (2pt);
\filldraw [gray] (1,2) circle (2pt);
\filldraw [gray] (2,2) circle (2pt);
\filldraw [gray] (3,2) circle (2pt);
\filldraw [gray] (4,2) circle (2pt);
\filldraw [gray] (1,3) circle (2pt);
\filldraw [gray] (2,3) circle (2pt);
\filldraw [gray] (3,3) circle (2pt);
\filldraw [gray] (4,3) circle (2pt);
\filldraw [gray] (2,3) circle (2pt);

\draw[very thick] (1,1) -- (1,1) -- (1,3);
\draw[very thick] (2,1) -- (2,1) -- (2,3);
\draw[very thick] (3,1) -- (3,1) -- (3,3);
\draw[very thick] (4,1) -- (4,1) -- (4,3);
\foreach \i in {3}
         {
           \draw[thick,dash dot,purple] (1,\i) .. controls (1.5,\i+.5) .. (2,\i);
           \draw[thick,dash dot,purple] (2,\i) .. controls (2.5,\i+.5) .. (3,\i);
           \draw[thick,dash dot,purple] (3,\i) .. controls (3.5,\i+.5) .. (4,\i);
           \draw[thick,dash dot,purple] (2,\i) .. controls (3,\i-.5) .. (4,\i);
         }
         
  \begin{pgfonlayer}{background}
    \filldraw [line width=4mm,black!3]
      (0.2,0.2)  rectangle (4.8,3.8);
  \end{pgfonlayer}
\end{tikzpicture}}
\hfill
\subcaptionbox{Diagram $\Gamma(\sigma,\pi)$ with order $\lambda^{5-6\alpha}=\lambda^{4-4\alpha}\lambda^{1-2\alpha}$.}[.49\textwidth]{
\begin{tikzpicture}[scale=0.9] 
\draw[black, thick] (0,0) rectangle (5,4);

\node[anchor=east,font=\small] at (0.8,3) {1};
\node[anchor=east,font=\small] at (0.8,2) {2};
\node[anchor=east,font=\small] at (0.8,1) {3};
\node[anchor=south,font=\small] at (1,0) {1};
\node[anchor=south,font=\small] at (2,0) {2};
\node[anchor=south,font=\small] at (3,0) {3};
\node[anchor=south,font=\small] at (4,0) {4};

\filldraw [gray] (1,1) circle (2pt);
\filldraw [gray] (2,1) circle (2pt);
\filldraw [gray] (3,1) circle (2pt);
\filldraw [gray] (4,1) circle (2pt);
\filldraw [gray] (1,2) circle (2pt);
\filldraw [gray] (2,2) circle (2pt);
\filldraw [gray] (3,2) circle (2pt);
\filldraw [gray] (4,2) circle (2pt);
\filldraw [gray] (1,3) circle (2pt);
\filldraw [gray] (2,3) circle (2pt);
\filldraw [gray] (3,3) circle (2pt);
\filldraw [gray] (4,3) circle (2pt);
\filldraw [gray] (2,3) circle (2pt);

\draw[very thick] (1,1) -- (1,1) -- (1,3);
\draw[very thick] (2,1) -- (2,2);
\draw[very thick] (3,1) -- (3,1) -- (3,3);
\draw[very thick] (4,1) -- (4,1) -- (4,3);

\draw[thick,dash dot,purple] (1,2) .. controls (1.5,2.5) .. (2,2);
\draw[thick,dash dot,purple] (2,2) .. controls (2.5,2.5) .. (3,2);
\draw[thick,dash dot,purple] (2,2) .. controls (3,1.5) .. (4,2);

\foreach \i in {3}
         {
           \draw[thick,dash dot,purple] (1,\i) .. controls (1.5,\i+.5) .. (2,\i);
           \draw[thick,dash dot,purple] (2,\i) .. controls (2.5,\i+.5) .. (3,\i);
           \draw[thick,dash dot,purple] (3,\i) .. controls (3.5,\i+.5) .. (4,\i);
           \draw[thick,dash dot,purple] (2,\i) .. controls (3,\i-.5) .. (4,\i);
         }
         
  \begin{pgfonlayer}{background}
    \filldraw [line width=4mm,black!3]
      (0.2,0.2)  rectangle (4.8,3.8);
  \end{pgfonlayer}
\end{tikzpicture}}
\caption{Splitting of a vertex with addition of three edges and $n=3$, $r=4$.}
\end{figure}

\vskip-0.3cm

\noindent

\smallskip

\begin{figure}[H]
\captionsetup[subfigure]{font=footnotesize}
\centering
\subcaptionbox{Diagram $\Gamma(\rho,\pi)$ with order $\lambda^{4-4\alpha}$.}[.49\textwidth]{\begin{tikzpicture}[scale=0.9] 
\draw[black, thick] (0,0) rectangle (5,4);

\node[anchor=east,font=\small] at (0.8,3) {1};
\node[anchor=east,font=\small] at (0.8,2) {2};
\node[anchor=east,font=\small] at (0.8,1) {3};
\node[anchor=south,font=\small] at (1,0) {1};
\node[anchor=south,font=\small] at (2,0) {2};
\node[anchor=south,font=\small] at (3,0) {3};
\node[anchor=south,font=\small] at (4,0) {4};

\filldraw [gray] (1,1) circle (2pt);
\filldraw [gray] (2,1) circle (2pt);
\filldraw [gray] (3,1) circle (2pt);
\filldraw [gray] (4,1) circle (2pt);
\filldraw [gray] (1,2) circle (2pt);
\filldraw [gray] (2,2) circle (2pt);
\filldraw [gray] (3,2) circle (2pt);
\filldraw [gray] (4,2) circle (2pt);
\filldraw [gray] (1,3) circle (2pt);
\filldraw [gray] (2,3) circle (2pt);
\filldraw [gray] (3,3) circle (2pt);
\filldraw [gray] (4,3) circle (2pt);
\filldraw [gray] (2,3) circle (2pt);

\draw[very thick] (1,1) -- (1,1) -- (1,3);
\draw[very thick] (2,1) -- (2,1) -- (2,3);
\draw[very thick] (3,1) -- (3,1) -- (3,3);
\draw[very thick] (4,1) -- (4,1) -- (4,3);
\foreach \i in {3}
         {
           \draw[thick,dash dot,purple] (1,\i) .. controls (1.5,\i+.5) .. (2,\i);
           \draw[thick,dash dot,purple] (2,\i) .. controls (2.5,\i+.5) .. (3,\i);
           \draw[thick,dash dot,purple] (3,\i) .. controls (3.5,\i+.5) .. (4,\i);
           \draw[thick,dash dot,purple] (2,\i) .. controls (3,\i-.5) .. (4,\i);
         }
         
  \begin{pgfonlayer}{background}
    \filldraw [line width=4mm,black!3]
      (0.2,0.2)  rectangle (4.8,3.8);
  \end{pgfonlayer}
\end{tikzpicture}}
\hfill
\subcaptionbox{Diagram $\Gamma(\sigma,\pi)$ with order $\lambda^{5-6\alpha}=\lambda^{4-4\alpha}\lambda^{1-2\alpha}$.}[.49\textwidth]{
\begin{tikzpicture}[scale=0.9] 
\draw[black, thick] (0,0) rectangle (5,4);

\node[anchor=east,font=\small] at (0.8,3) {1};
\node[anchor=east,font=\small] at (0.8,2) {2};
\node[anchor=east,font=\small] at (0.8,1) {3};
\node[anchor=south,font=\small] at (1,0) {1};
\node[anchor=south,font=\small] at (2,0) {2};
\node[anchor=south,font=\small] at (3,0) {3};
\node[anchor=south,font=\small] at (4,0) {4};

\filldraw [gray] (1,1) circle (2pt);
\filldraw [gray] (2,1) circle (2pt);
\filldraw [gray] (3,1) circle (2pt);
\filldraw [gray] (4,1) circle (2pt);
\filldraw [gray] (1,2) circle (2pt);
\filldraw [gray] (2,2) circle (2pt);
\filldraw [gray] (3,2) circle (2pt);
\filldraw [gray] (4,2) circle (2pt);
\filldraw [gray] (1,3) circle (2pt);
\filldraw [gray] (2,3) circle (2pt);
\filldraw [gray] (3,3) circle (2pt);
\filldraw [gray] (4,3) circle (2pt);
\filldraw [gray] (2,3) circle (2pt);

\draw[very thick] (1,1) -- (1,1) -- (1,3);
\draw[very thick] (2,1) -- (2,1) -- (2,3);
\draw[very thick] (3,1) -- (3,2);
\draw[very thick] (4,1) -- (4,1) -- (4,3);

\draw[thick,dash dot,purple] (2,2) .. controls (2.5,2.5) .. (3,2);
\draw[thick,dash dot,purple] (3,2) .. controls (3.5,2.5) .. (4,2);

\foreach \i in {3}
         {
           \draw[thick,dash dot,purple] (1,\i) .. controls (1.5,\i+.5) .. (2,\i);
           \draw[thick,dash dot,purple] (2,\i) .. controls (2.5,\i+.5) .. (3,\i);
           \draw[thick,dash dot,purple] (3,\i) .. controls (3.5,\i+.5) .. (4,\i);
           \draw[thick,dash dot,purple] (2,\i) .. controls (3,\i-.5) .. (4,\i);
         }
         
  \begin{pgfonlayer}{background}
    \filldraw [line width=4mm,black!3]
      (0.2,0.2)  rectangle (4.8,3.8);
  \end{pgfonlayer}
\end{tikzpicture}}
\caption{Splitting of a vertex with addition of two edges and $n=3$, $r=4$.}
\end{figure}

\vskip-0.3cm

\noindent

\smallskip

\begin{figure}[H]
\captionsetup[subfigure]{font=footnotesize}
\centering
\subcaptionbox{Diagram $\Gamma(\rho,\pi)$ with order $\lambda^{4-4\alpha}$.}[.49\textwidth]{
\begin{tikzpicture}[scale=0.9] 
\draw[black, thick] (0,0) rectangle (5,4);

\node[anchor=east,font=\small] at (0.8,3) {1};
\node[anchor=east,font=\small] at (0.8,2) {2};
\node[anchor=east,font=\small] at (0.8,1) {3};
\node[anchor=south,font=\small] at (1,0) {1};
\node[anchor=south,font=\small] at (2,0) {2};
\node[anchor=south,font=\small] at (3,0) {3};
\node[anchor=south,font=\small] at (4,0) {4};

\filldraw [gray] (1,1) circle (2pt);
\filldraw [gray] (2,1) circle (2pt);
\filldraw [gray] (3,1) circle (2pt);
\filldraw [gray] (4,1) circle (2pt);
\filldraw [gray] (1,2) circle (2pt);
\filldraw [gray] (2,2) circle (2pt);
\filldraw [gray] (3,2) circle (2pt);
\filldraw [gray] (4,2) circle (2pt);
\filldraw [gray] (1,3) circle (2pt);
\filldraw [gray] (2,3) circle (2pt);
\filldraw [gray] (3,3) circle (2pt);
\filldraw [gray] (4,3) circle (2pt);
\filldraw [gray] (2,3) circle (2pt);

\draw[very thick] (1,1) -- (1,1) -- (1,3);
\draw[very thick] (2,1) -- (2,1) -- (2,3);
\draw[very thick] (3,1) -- (3,1) -- (3,3);
\draw[very thick] (4,1) -- (4,1) -- (4,3);
\foreach \i in {3}
         {
           \draw[thick,dash dot,purple] (1,\i) .. controls (1.5,\i+.5) .. (2,\i);
           \draw[thick,dash dot,purple] (2,\i) .. controls (2.5,\i+.5) .. (3,\i);
           \draw[thick,dash dot,purple] (3,\i) .. controls (3.5,\i+.5) .. (4,\i);
           \draw[thick,dash dot,purple] (2,\i) .. controls (3,\i-.5) .. (4,\i);
         }
         
  \begin{pgfonlayer}{background}
    \filldraw [line width=4mm,black!3]
      (0.2,0.2)  rectangle (4.8,3.8);
  \end{pgfonlayer}
\end{tikzpicture}}
\hfill
\subcaptionbox{Diagram $\Gamma(\sigma,\pi)$ with order $\lambda^{5-6\alpha}=\lambda^{4-4\alpha}\lambda^{1-2\alpha}$.}[.49\textwidth]{
\begin{tikzpicture}[scale=0.9] 
\draw[black, thick] (0,0) rectangle (5,4);

\node[anchor=east,font=\small] at (0.8,3) {1};
\node[anchor=east,font=\small] at (0.8,2) {2};
\node[anchor=east,font=\small] at (0.8,1) {3};
\node[anchor=south,font=\small] at (1,0) {1};
\node[anchor=south,font=\small] at (2,0) {2};
\node[anchor=south,font=\small] at (3,0) {3};
\node[anchor=south,font=\small] at (4,0) {4};

\filldraw [gray] (1,1) circle (2pt);
\filldraw [gray] (2,1) circle (2pt);
\filldraw [gray] (3,1) circle (2pt);
\filldraw [gray] (4,1) circle (2pt);
\filldraw [gray] (1,2) circle (2pt);
\filldraw [gray] (2,2) circle (2pt);
\filldraw [gray] (3,2) circle (2pt);
\filldraw [gray] (4,2) circle (2pt);
\filldraw [gray] (1,3) circle (2pt);
\filldraw [gray] (2,3) circle (2pt);
\filldraw [gray] (3,3) circle (2pt);
\filldraw [gray] (4,3) circle (2pt);
\filldraw [gray] (2,3) circle (2pt);

\draw[very thick] (1,1) -- (1,1) -- (1,3);
\draw[very thick] (2,1) -- (2,1) -- (2,3);
\draw[very thick] (3,1) -- (3,1) -- (3,3);
\draw[very thick] (4,1) -- (4,2);

\draw[thick,dash dot,purple] (3,2) .. controls (3.5,2.5) .. (4,2);
\draw[thick,dash dot,purple] (2,2) .. controls (3,1.5) .. (4,2);

\foreach \i in {3}
         {
           \draw[thick,dash dot,purple] (1,\i) .. controls (1.5,\i+.5) .. (2,\i);
           \draw[thick,dash dot,purple] (2,\i) .. controls (2.5,\i+.5) .. (3,\i);
           \draw[thick,dash dot,purple] (3,\i) .. controls (3.5,\i+.5) .. (4,\i);
           \draw[thick,dash dot,purple] (2,\i) .. controls (3,\i-.5) .. (4,\i);
         }
         
  \begin{pgfonlayer}{background}
    \filldraw [line width=4mm,black!3]
      (0.2,0.2)  rectangle (4.8,3.8);
  \end{pgfonlayer}
\end{tikzpicture}}
\caption{Splitting of a vertex with addition of two edges and $n=3$, $r=4$.}
\label{fig:diagram3-1-3}
\end{figure}

\vskip-0.3cm

\noindent
 When $G$ is a triangle with $n=2$ and $r=3$, 
 the above procedure can be reversed 
 by first merging a vertex and then gluing edges, see Figure~\ref{fig:diagram3-23},
 which results into ``overlapping'' all copies of the graph $G$. 
  
\begin{figure}[H]
\captionsetup[subfigure]{font=footnotesize}
\centering
\subcaptionbox{Merging one vertex.}[.3\textwidth]{
\begin{tikzpicture}[scale=0.9] 
\draw[black, thick] (0.2,0.2) rectangle (3.8,4.3);

\filldraw [gray] (1,1) circle (2pt);
\filldraw [gray] (2,1) circle (2pt);
\filldraw [gray] (3,1) circle (2pt);
\filldraw [gray] (1,2) circle (2pt);
\filldraw [gray] (2,2) circle (2pt);
\filldraw [gray] (3,2) circle (2pt);

\draw[very thick] (1,1) -- (1,2);
\foreach \i in {1,2}
         {
           \draw[thick,dash dot,purple] (1,\i) .. controls (1.5,\i+.5) .. (2,\i);
           \draw[thick,dash dot,purple] (2,\i) .. controls (2.5,\i+.5) .. (3,\i);
           \draw[thick,dash dot,purple] (1,\i) .. controls (2,\i-.5) .. (3,\i);
         }
\filldraw [black] (1.1,3.8) circle (1.5pt);         
\filldraw [black] (1.1,2.8) circle (1.5pt);   
\filldraw [black] (1.97,3.3) circle (1.5pt); 
\filldraw [black] (2.84,3.8) circle (1.5pt);  
\filldraw [black] (2.84,2.8) circle (1.5pt);  
\draw[thick,blue] (1.1,3.8) -- (2.84,2.8);
\draw[thick,blue] (1.1,2.8) -- (1.1,3.8);
\draw[thick,blue] (1.1,2.8) -- (2.84,3.8);
\draw[thick,blue] (2.84,2.8) -- (2.84,3.8);
         
  \begin{pgfonlayer}{background}
    \filldraw [line width=4mm,black!3]
      (0.4,0.4)  rectangle (3.6,4.1);
  \end{pgfonlayer}
\end{tikzpicture}}
\hfill
\subcaptionbox{Gluing one edge.}[.3\textwidth]{
\begin{tikzpicture}[scale=0.9] 
\draw[black, thick] (0.2,0.2) rectangle (3.8,4.3);

\filldraw [gray] (1,1) circle (2pt);
\filldraw [gray] (2,1) circle (2pt);
\filldraw [gray] (3,1) circle (2pt);
\filldraw [gray] (1,2) circle (2pt);
\filldraw [gray] (2,2) circle (2pt);
\filldraw [gray] (3,2) circle (2pt);

\draw[very thick] (1,1) -- (1,2);
\draw[very thick] (3,1) -- (3,2);
\foreach \i in {1,2}
         {
           \draw[thick,dash dot,purple] (1,\i) .. controls (1.5,\i+.5) .. (2,\i);
           \draw[thick,dash dot,purple] (2,\i) .. controls (2.5,\i+.5) .. (3,\i);
         }
          \draw[thick,dash dot,purple] (1,2) .. controls (2,1.5) .. (3,2);
\filldraw [black] (2,3.8) circle (1.5pt);         
\filldraw [black] (2,2.8) circle (1.5pt);   
\filldraw [black] (1.13,3.3) circle (1.5pt);   
\filldraw [black] (2.87,3.3) circle (1.5pt);  
 
\draw[thick,blue] (2,3.8) -- (2,2.8);
\draw[thick,blue] (2,3.8) -- (1.13,3.3);
\draw[thick,blue] (2,2.8) -- (1.13,3.3);
\draw[thick,blue] (2,3.8) -- (2.87,3.3);
\draw[thick,blue] (2,2.8) -- (2.87,3.3);
                 
  \begin{pgfonlayer}{background}
    \filldraw [line width=4mm,black!3]
      (0.4,0.4)  rectangle (3.6,4.1);
  \end{pgfonlayer}
\end{tikzpicture}}
\hfill
\subcaptionbox{Gluing three edges.}[.3\textwidth]{
\begin{tikzpicture}[scale=0.9] 
\draw[black, thick] (0.2,0.2) rectangle (3.8,4.3);

\filldraw [gray] (1,1) circle (2pt);
\filldraw [gray] (2,1) circle (2pt);
\filldraw [gray] (3,1) circle (2pt);
\filldraw [gray] (1,2) circle (2pt);
\filldraw [gray] (2,2) circle (2pt);
\filldraw [gray] (3,2) circle (2pt);
\draw[very thick] (1,1) -- (1,2);
\draw[very thick] (2,1) -- (2,2);
\draw[very thick] (3,1) -- (3,2);
\foreach \i in {2}
         {
           \draw[thick,dash dot,purple] (1,\i) .. controls (1.5,\i+.5) .. (2,\i);
           \draw[thick,dash dot,purple] (2,\i) .. controls (2.5,\i+.5) .. (3,\i);
           \draw[thick,dash dot,purple] (1,\i) .. controls (2,\i-.5) .. (3,\i);
         }

\filldraw [black] (2,3.8) circle (1.5pt);         
\filldraw [black] (1.5,2.93) circle (1.5pt);   
\filldraw [black] (2.5,2.93) circle (1.5pt);   

\draw[thick,blue] (2,3.8) -- (1.5,2.93);
\draw[thick,blue] (2,3.8) -- (2.5,2.93);
\draw[thick,blue] (2.5,2.93) -- (1.5,2.93);
                  
  \begin{pgfonlayer}{background}
    \filldraw [line width=4mm,black!3]
      (0.4,0.4)  rectangle (3.6,4.1);
  \end{pgfonlayer}
\end{tikzpicture}}
\caption{Diagram patterns with $G$ a triangle and $n=2$, $r=3$.}
\label{fig:diagram3-23}
\end{figure}

\vskip-0.3cm

\noindent
 As in part~$(b)$ above, we lower
 bound $\kappa_n(N_G)$ using a single partition,
 and we upper bound using the total count of
 connected non-flat partitions using Lemma~\ref{fjkldsf-l}-$b)$
 to obtain \eqref{cumulant-rhoph1}.
 
\noindent
   $c)$ is a direct consequence of part~$(b)$ above.
\end{Proof}

\begin{corollary}[Sparse regime] 
\label{th6.4-c}
    Let $G$ be a connected graph with $|V(G)|=r$ vertices, $r\geq 2$,
  satisfying Assumption~\ref{a61} for $n=2$ in the sparse regime \eqref{fjnldsf-2}. 
  \begin{enumerate}[a)] 
  \item 
    If $G$ is a tree, i.e. $|E(G)| = r-1$, we have the normalized 
    cumulant bounds 
 \begin{equation}
   \label{jfkla} 
   \big|\kappa_n(\widetilde{N}_G)\big|
  \leq 
   (K_3)^n 
 n!^r
      \lambda^{
       - (
\alpha       -(\alpha - 1)r 
       ) ( n/2-1 ) },
   \qquad \lambda \geq 1,
   \quad n \geq 2, 
\end{equation} 
 where $K_3:=\max ( (K_2)^r , 1) /(K_1)^{r/2}$.
\item
 If $G$ is not a tree, i.e. $|E(G)|\geq r$,
 we have the normalized cumulant bounds 
 \begin{equation}
   \label{fjklds34} 
    \big|\kappa_n(\widetilde{N}_G)\big|
  \leq 
    n!^r
      (K_3)^n
        \lambda^{(\alpha |E(G)|-r)(n/2-1)}  , 
   \qquad \lambda \geq 1,
  \quad n \geq 2, 
\end{equation} 
 for some $K_3 > 0$. 
\item
  If $G$ is a cycle, i.e. $|E(G)| = r$,
  we have the normalized cumulant bounds 
 \begin{equation}
   \label{fjklds34-2} 
    \big|\kappa_n(\widetilde{N}_G)\big|
  \leq 
  n!^r
   (K_3)^n
    \lambda^{ (\alpha -1 )(n/2-1)r}, 
     \qquad \lambda \geq 1,
 \quad n \geq 2, 
\end{equation} 
 for some $K_3>0$.
  \end{enumerate}
\end{corollary}
\begin{Proof}
  We note that the upper bound in \eqref{equiv-1}
  does not require Assumption~\ref{a61}. 
  Regarding \eqref{jfkla}, we have  
\begin{eqnarray*} 
  \big|\kappa_n(\widetilde{N}_G)\big|
& \leq &  
  \frac{n!^r 
    (K_2)^r}{
   (
    (K_1)^r
        \lambda^{
\alpha      -(\alpha - 1)r 
        }
        )^{n/2}}
  \lambda^{
 \alpha          -(\alpha - 1)r 
  }
    \\
    & = & 
   \frac{(K_2)^r}{(K_1)^{nr/2}}
   n!^r
        \lambda^{-(
 \alpha      -(\alpha - 1)r 
       ) ( n/2-1)
   }, \qquad n \geq 2.
\end{eqnarray*} 
   \noindent
 Regarding \eqref{fjklds34}, we have  
$$ 
   \big|\kappa_n(\widetilde{N}_G)\big|
 \leq 
   \frac{
  n!^r
    (K_2)^r
  \lambda^{r-\alpha |E(G)|}
   }{
            (
     (K_1)^r
     \lambda^{r-\alpha |E(G)|}
     )^{n/2}}
=
   n!^r
   \frac{
       (K_2)^r }{(K_1)^{nr/2}}
      \lambda^{ - (r-\alpha |E(G)|)(n/2 - 1)}  , 
  \quad n \geq 2. 
$$ 
 Finally, the bound 
 \eqref{fjklds34-2} is a direct consequence of \eqref{fjklds34}. 
\end{Proof}
\section{Asymptotic normality of subgraph counts}
\label{s6-1}
\noindent
 In this section, we let $H(x,y)$ be a connection function  
  satisfying Assumption~\ref{a61},
  and consider the
  random-connection model $G_{H_\lambda} (\Xi)$
  where $H_\lambda(x,y):= c_\lambda H(x,y)$,
  $\lambda >0$.
\begin{corollary}[Dilute regime]
  \label{c01}
  Let $G$ be a connected graph with $|V(G)|=r$ vertices, $r\geq 2$,
  satisfying Assumption~\ref{a61}. 
  In the dilute regime \eqref{fjnldsf},
  the normalized subgraph count $\widetilde{N}_G$ in 
  $G_{H_\lambda} (\Xi)$ satisfies
  the Kolmogorov distance bound 
\begin{equation}
  \label{fjkld13}
  \sup_{x\in \real}
\big| \P \big( \widetilde{N}_G \leq x \big) - \Phi(x) \big| \leq
C \lambda^{ - 1/(4r - 2)},
\end{equation}
with rate
{exponent}
$1/(4r -2)$ as $\lambda$ tends to infinity, where $C>0$ depends only on $H$ and $G$.
\end{corollary}
\begin{Proof}
   In the dilute regime, the cumulant bound 
 \eqref{Statuleviciuscond} 
 shows that
 the centered and normalized subgraph count
 $\widetilde{N}_G$ 
 satisfies the {Statulevi\v{c}ius condition}
 \eqref{Statuleviciuscond2} in {the}
 appendix, see \cite{rudzkis,doering},
 with $\gamma := r-1$. 
 We conclude by applying Corollary~\ref{t1-c}
 and Lemma~\ref{l1}-$i)$ with $\gamma :=r-1$
 and $\Delta_\lambda:=\sqrt{K \lambda}$. 
\end{Proof}
In the sparse regime we have the following result, 
 in which \eqref{fjkldsa1} is consistent with
 \eqref{fjkld13} when $\alpha = 1$.
\begin{corollary}[Sparse regime]
  \label{c01-2}
  Let $G$ be a tree with $|V(G)| = r \geq 2$ vertices,
  satisfying Assumption~\ref{a61}. 
  In the sparse regime \eqref{fjnldsf-2}
  with $\alpha \in [1, r/(r-1) )$,
  the normalized subgraph count $\widetilde{N}_G$ in 
  $G_{H_\lambda} (\Xi)$ satisfies
  the Kolmogorov distance bound 
\begin{equation}
\label{fjkldsa1} 
\sup_{x\in \real}
\big| \P \big( \widetilde{N}_G \leq x \big) - \Phi(x) \big| \leq
C \lambda^{ - (
 \alpha   -(\alpha - 1)r 
    ) / ( 4r - 2) }, 
\end{equation} 
 as $\lambda$ tends to infinity, where $C>0$ depends only on $H$ and $G$.
\end{corollary}
\begin{Proof}
  This is a consequence of 
  Corollary~\ref{th6.4-c}-$a)$
  and Lemma~\ref{l1}-$i)$ in {the} appendix,
  with $\gamma :=r-1$ and 
  $\Delta_\lambda:=
  (K \lambda )^{ - (
 \alpha   -(\alpha - 1)r 
    )/2}$
  and $\alpha \in [1, r/(r-1) )$. 
\end{Proof} 
  We note that  up to division by $2r - 1$,
   the rate in \eqref{fjkldsa1} is consistent  
   with the rate
   {exponent}
   $(
 \alpha    -(\alpha - 1)r 
     ) / 2$ obtained for the counting of trees in the
 Erd{\H o}s-R\'enyi graph, cf. Corollary~4.10 of \cite{PS2}. 
 In addition, since 
$
(\alpha |E(G)|-r)(n/2-1)
\geq
(\alpha -1)(n/2-1)r \geq 0$,
no significant Kolmogorov bounds
are derived from \eqref{fjklds34}
and \eqref{fjklds34-2} 
for cycle and other non-tree graphs 
in the sparse regime,
which is consistent with
Corollaries~4.8-4.9 of \cite{PS2}. 

\medskip

 Taking $\Delta_\lambda = \sqrt{K \lambda}$, 
 by Lemma~\ref{l1}-$ii)$ in {the} appendix, see Theorem~1.1 of \cite{doring},
we have the following result. 
\begin{corollary}
  \label{c01-2-0}
  Let $G$ be a connected graph with $|V(G)| = r \geq 2$ vertices,
  satisfying Assumption~\ref{a61}. 
  The normalized subgraph count $\widetilde{N}_G$
    satisfies a moderate deviation principle
  in the dilute regime of Corollary~\ref{c01}, 
  with speed $a_\lambda^2 = o( \lambda^{1/(2r - 1)} )$ and rate function $x^2/2$.
\end{corollary}
In addition,
by Lemma~\ref{l1}-$iii)$ in {the} appendix,
see the corollary of \cite[Lemma~2.4]{saulis},
{there exists a constant $K>0$
such that for any sufficiently large $\lambda$ 
}
we have the concentration inequality 
\begin{equation}
\label{concentrationineq}
  \P \big( \big| \widetilde{N}_G \big|
  \geq x)\le2\exp\left(-
  \frac{1}{4} 
  {
    \min \left(
  \frac{x^2}{2^r}
  ,
  \big( x \sqrt{K \lambda }\big)^{1/r}
  \right)
  }
  \right), 
  \quad
  x\ge0, 
\end{equation} 
in agreement with the rate in Theorem~1.1 of \cite{bachmann},
which is stated for subgraph counts in random geometric
graphs. 
\section{Subgraph counts in random geometric graphs}
\label{rgg}
\noindent
In this section, we consider subgraph counts in
the (Poisson) random geometric graph model.
Assume that $\mu$ is the Lebesgue measure,
and that the intensity measure $\Lambda$ takes the form 
$$\Lambda (\mathrm{d}x) := {\bf 1}_A (x) \mu(\mathrm{d}x),\qquad\lambda>0,$$
  where $A$ is a Borel subset of $\real^d$
  such that $\mu (A)<\infty$. 
\begin{definition}
  For every $\lambda >0$,
  let $G_{H_\lambda} (\Xi)$
  denote the random-connection model
  with connection function 
  $$
  H_\lambda (x,y):=\bone_{\{\|x-y\|\leq R_\lambda\}},
  \qquad x,y\in\R^d,
  $$
  for some function $R_\lambda>0$ of $\lambda$.
 We consider the following regimes. 
\begin{itemize}
\item Dense regime: we have 
$$ 
    \lim_{\lambda\to \infty} \lambda R_\lambda^d \in (0,\infty]. 
$$
\item Sparse regime: we have 
$$ 
 \lim_{\lambda\to \infty} \lambda R_\lambda^d = 0
\ \mathrm{and} \
\lim_{\lambda\to \infty}
\lambda \big( \lambda R_\lambda^d \big)^{r-1} = \infty.
$$
\end{itemize} 
\end{definition}
When $\lim_{\lambda\to \infty} \lambda R_\lambda^d = c \in (0,\infty )$, 
we also say that we are in the thermodynamic regime.
The following result extends Proposition~3.2 of \cite{lachiezerey2}
from second order cumulants to cumulants of any order. 
\begin{prop}
    \label{thm8}
 Let $G$ be a connected graph with $|V(G)|=r$ vertices, $r\geq 2$.
 In the random geometric graph model we have the following cumulant bounds. 
 \begin{enumerate}[a)]
 \item (Dense regime). 
   We have
   \begin{equation}
     \label{b1}
     K_1 (n-1)! \lambda^{1+(r-1)n} ( R_\lambda^d)^{(r-1)n}
     \leq \kappa_n(N_G) \leq K_2 n!^r r!^{n-1}
 \lambda^{1+(r-1)n} ( R_\lambda^d)^{(r-1)n}, 
 \quad \lambda \geq 1,
   \end{equation}
for some constants $K_1, K_2 > 0$ independent
of $\lambda , n\geq 1$. 
  \item (Thermodynamic regime). 
   We have
   \begin{equation}
\label{b2}
K_1 \lambda \leq \kappa_n(N_G)\leq K_2 n!^r r!^{n-1}\lambda,
  \quad \lambda \geq 1, 
\end{equation}
 for some constants $K_1, K_2 > 0$ independent
of $\lambda , n\geq 1$. 
  \item (Sparse regime). 
   We have
   \begin{equation}
 \label{b3}
         K_1 \lambda^r ( R_\lambda^d)^{r-1} 
     \leq \kappa_n(N_G)\leq K_2 n!^r r!^{n-1} \lambda^r (R_\lambda^d)^{r-1},
     \quad \lambda \geq 1,
   \end{equation}
 for some constants $K_1, K_2 > 0$ independent
of $\lambda , n\geq 1$. 
 \end{enumerate}
\end{prop}
\begin{Proof}
 By Proposition~\ref{lma-diagram1},
 letting $\widebar{\rho}$ denote a spanning tree contained in $\rho$,
  we have 
\begin{eqnarray}
  \nonumber
  \kappa_n(N_G)&=&\sum_{\rho \in \Pi_{\widehat{1}} ( [n] \times [r])
    \atop
    {\rho \wedge \pi = \widehat{0}}
  }
  \lambda^{|\rho|}\int_{A^{|\rho |}}
  \left(
\prod_{\{i,j\} \in E(\rho_G)}
\bone_{\{\|x_i-x_j\|\leq R_\lambda\}}
\right)
\mathrm{d}x_1\cdots\mathrm{d}x_{|\rho|}
\\
\nonumber
&\leq &\sum_{\rho \in \Pi_{\widehat{1}} ( [n] \times [r])
    \atop
    {\rho \wedge \pi = \widehat{0}}
  }
  \lambda^{|\rho|}\int_{A^{|\rho |}}
  \left(
\prod_{\{i,j\} \in E(\widebar{\rho}_G)}
\bone_{\{\|x_i-x_j\|\leq R_\lambda\}}
\right)
\mathrm{d}x_1\cdots\mathrm{d}x_{|\rho|}
\\
\nonumber
&=&\sum_{\rho \in \Pi_{\widehat{1}} ( [n] \times [r])
    \atop
    {\rho \wedge \pi = \widehat{0}}
  }
  \lambda^{|\rho|} \mu (A) (v_d R_\lambda^d)^{|\rho|-1}. 
\end{eqnarray}
 
  \noindent
$a)$ 
  In the dense regime
  with $\lim_{\lambda\to \infty} \lambda R_\lambda^d = \infty$,
  the dominating asymptotic order 
  $\lambda(\lambda R_\lambda^d)^{(r-1)n}$ of $\kappa_n(N_G)$ is
  achieved when $|\rho|=1+(r-1)n$,
  which yields the upper bound in \eqref{b1}. 

\noindent
$b)$ 
    In the thermodynamic regime
    with $\lim_{\lambda\to \infty} \lambda R_\lambda^d = c>0$,
    the dominating asymptotic order of $\kappa_n(N_G)$ is $\lambda$,
    which yields the upper bound in \eqref{b2}. 

  \noindent
$c)$ 
    In the sparse regime
    with $\lim_{\lambda\to \infty} \lambda R_\lambda^d = 0$ and
    $\lim_{\lambda\to \infty} \lambda(\lambda R_\lambda^d)^{r-1} = \infty$,
    the dominating asymptotic order $\lambda(\lambda R_\lambda^d)^{r-1}$
    of $\kappa_n(N_G)$ is achieved when $|\rho|=r$, 
    which yields the upper bound in \eqref{b3}. 

    \medskip 

 In addition, the kernel 
 $H_\lambda (x,y)=\bone_{\{\|x-y\|\leq R_\lambda\}}$
 satisfies Assumption~\ref{a61} for all $n\geq 1$,
 with $C_{H_\lambda} = v_d (R_\lambda /2)^d$
 in the framework of
 {above increasing intensity example},
 which similarly yields the lower bounds in 
 \eqref{b1}-\eqref{b3}.
\end{Proof}
 The next result is a direct consequence of Proposition~\ref{thm8}. 
\begin{corollary}
    \label{thm8-0}
 Let $G$ be a connected graph with $|V(G)|=r$ vertices, $r\geq 2$.
 In the random geometric graph model we have the following
 normalized cumulant bounds.  
 \begin{enumerate}[a)]
 \item (Dense and thermodynamic regime). 
   We have
   \begin{equation}
     \kappa_n\big(\widetilde{N}_G\big) \leq n!^r ( K \lambda )^{-(n/2-1)}, 
     \quad \lambda \geq 1,
   \end{equation}
for some $K>0$ constant independent of $\lambda , n\geq 1$. 
  \item (Sparse regime). 
   We have
   \begin{equation}
     \kappa_n\big(\widetilde{N}_G\big)\leq K \big(
      \lambda^r R_\lambda^{(r-1)d} \big)^{-(n/2-1)},
      \quad \lambda \geq 1,
   \end{equation}
 for some constants $K_1, K_2 > 0$ independent of $\lambda , n\geq 1$. 
 \end{enumerate}
\end{corollary}
 The following result then follows from Corollary~\ref{thm8-0}. 
\begin{corollary}
  \label{jdkj10}
  Let $G$ be a connected graph with $|V(G)|=r$ vertices, $r\geq 2$,
  in the random geometric graph model.
  \begin{enumerate}[i)]
    \item 
    Dense \!\! \slash \ thermodynamic regimes. 
  If $\lim_{\lambda\to\infty} (\lambda R_\lambda^d) \in (0,\infty]$, 
 we have
\begin{equation}
\nonumber
  \sup_{x\in \real}
\big| \P \big( \widetilde{N}_G \leq x \big) - \Phi(x) \big| \leq
C
\lambda^{-1/ ( 4r - 2) }.
\end{equation}
\item
  Sparse regime.
  If $\lim_{\lambda\to\infty} \lambda R_\lambda^d = 0$ and
  $\lim_{\lambda\to\infty} \lambda(\lambda R_\lambda^d)^{r-1} = \infty$,
  we have 
\begin{equation}
\sup_{x\in \real}
\big| \P \big( \widetilde{N}_G \leq x \big) - \Phi(x) \big| \leq
C
\big(\lambda^r R_\lambda^{(r-1)d}\big)^{-1/ ( 4r - 2) }.
\end{equation}
\end{enumerate}
\end{corollary}
\begin{Proof}
  In both cases $(i)$ and $(ii)$ we apply Corollary~\ref{thm8-0} 
  and Lemma~\ref{l1}-$i)$ with $\gamma :=r-1$,
  by taking 
  $\Delta_\lambda:=
  \sqrt{\lambda}$
  in the dense and  thermodynamic regimes,
  and 
  $\Delta_\lambda:= \lambda^{r/2} R_\lambda^{(r-1)d/2}$
  in the sparse regime.
\end{Proof} 
We note that Berry-Esseen convergence rates have been obtained
for certain random functionals in the random geometric graph
model, including total edge lengths in \cite[Corollary~4.3]{schulte},
clique counts using Poisson $U$-statistics in \cite[Theorem~4.1]{reitzner}  
and using stabilizing functionals in \cite[Theorem~3.15]{lachiezerey4},
and $k$-hop counts in the one-dimensional unit disk model in
\cite[Proposition~8.1]{privaultkhops}. 

\subsubsection*{Moderate deviation and concentration inequalities} 
\noindent 
Letting $\gamma =r-1$, 
In the dense and thermodynamic regimes of 
Corollary~\ref{jdkj10} with $\Delta_\lambda=
\sqrt{\lambda}$, $\widetilde{N}_G/a_\lambda$
    satisfies a moderate deviation principle
with rate function $x^2/2$
and speed $a_\lambda^2 = o( \lambda^{1/(2r - 1)} )$
in the setting of Lemma~\ref{l1}-$ii)$ in {the} appendix,
and the concentration inequality 
\eqref{concentrationineq} holds by Lemma~\ref{l1}-$iii)$. 

\appendix

\section{Appendix}
\noindent 
The following results are summarized from the ``main lemmas'' in Chapter~2 of \cite{saulis} and \cite{doring}, and are tailored for our applications
to the random-connection model. 
\begin{lemma}
  \label{l1}
  Let $(X_\lambda)_{\lambda \geq 1}$ be a family of random variables with mean zero and unit variance for all $\lambda>0$. Suppose that for all $\lambda \geq 1$, all moments of the random variable $X_\lambda$ exist and that 
  the cumulants of $X_\lambda$ satisfy  
  \begin{equation}
    \label{Statuleviciuscond2}
    |\kappa_n (X_\lambda)|\le\frac{(n!)^{1+\gamma}}{(\Delta_\lambda)^{n-2}},
    \qquad
 n\ge3, 
\end{equation}
  where $\gamma\ge0$ is a constant not depending on $\lambda$, while $\Delta_\lambda\in(0,\infty)$ may depend on $\lambda$.
   Then, the following assertions hold.
\begin{enumerate}[i)]
\item (Kolmogorov bound,
 \cite[Corollary~2.1]{saulis} and \cite[Theorem~2.4]{doering})
  One has
\begin{equation}
\sup_{x\in\R}|\IP(X_\lambda\leq x)-\Phi(x)|\leq \frac{C}{(\Delta_\lambda)^{1/(1+2\gamma)}},
\end{equation}
for some constant $C>0$ depending only on $\gamma$.
\item (Moderate deviation principle,
  \cite[Theorem~1.1]{doring}
    {
    and \cite[Theorem~3.1]{doering}}).
  Let $( a_\lambda )_{\lambda > 0}$ be a sequence of real numbers tending to infinity, and such that 
  $$
  \lim_{\lambda \to \infty}
  \frac{a_\lambda}{(\Delta_\lambda)^{1/(1+2\gamma)}}
  = 0.
  $$
    Then, $ (a_\lambda^{-1}X_\lambda)_{\lambda >0}$ satisfies a moderate deviation principle with speed $a_\lambda^2$ and rate function $x^2/2$. 
\item (Concentration inequality,
  corollary of \cite[Lemma~2.4]{saulis}
  {
  and \cite[Theorem~2.5]{doering}}).
  For any sufficiently large $\lambda$, 
\begin{equation}
  \IP(|X_\lambda|\geq x)\le2\exp\left(-\frac14\min\left( \frac{x^2}{2^{
      {1+\gamma}}},(x\Delta_\lambda)^{1/(1+\gamma)}\right) \right),
  \quad
  x\ge0. 
\end{equation} 
\end{enumerate}
\end{lemma}

\subsubsection*{Acknowledgement}
\noindent
We thank M. Schulte and C. Th\"ale for a correction to an earlier
version of Lemma~\ref{fjkldsf-l}. 

\footnotesize

\def\cprime{$'$} \def\polhk#1{\setbox0=\hbox{#1}{\ooalign{\hidewidth
  \lower1.5ex\hbox{`}\hidewidth\crcr\unhbox0}}}
  \def\polhk#1{\setbox0=\hbox{#1}{\ooalign{\hidewidth
  \lower1.5ex\hbox{`}\hidewidth\crcr\unhbox0}}} \def\cprime{$'$}


\begin{thebibliography}{BRSW17}

\bibitem[BKR89]{BKR}
A.D. Barbour, M.~Karo{\'n}ski, and A.~Ruci{\'n}ski.
\newblock A central limit theorem for decomposable random variables with
  applications to random graphs.
\newblock {\em J. Combin. Theory Ser. B}, 47(2):125--145, 1989.

\bibitem[BOR85]{eabender}
E.A. Bender, A.M. Odlyzko, and L.B. Richmond.
\newblock The asymptotic number of irreducible partitions.
\newblock {\em European J. Combin.}, 6(1):1--6, 1985.

\bibitem[BR12]{balakrishnan}
R.~Balakrishnan and K.~Ranganathan.
\newblock {\em A Textbook of Graph Theory}.
\newblock Universitext. Springer, New York, second edition, 2012.

\bibitem[BR18]{bachmann}
S.~Bachmann and M.~Reitzner.
\newblock Concentration for {P}oisson {$U$}-statistics: {S}ubgraph counts in
  random geometric graphs.
\newblock {\em Stochastic Process. Appl.}, 128:3327--3352, 2018.

\bibitem[BRSW17]{bogdan}
K.~Bogdan, J.~Rosi\'{n}ski, G.~Serafin, and L.~Wojciechowski.
\newblock L\'{e}vy systems and moment formulas for mixed {P}oisson integrals.
\newblock In {\em Stochastic analysis and related topics}, volume~72 of {\em
  Progr. Probab.}, pages 139--164. Birkh{\"{a}}user/Springer, Cham, 2017.

\bibitem[CT22]{can2022}
V.~H. Can and K.~D. Trinh.
\newblock Random connection models in the thermodynamic regime: central limit
  theorems for add-one cost stabilizing functionals.
\newblock {\em Electron. J. Probab.}, 27:1--40, 2022.

\bibitem[DE13]{doring}
H.~D{\"{o}}ring and P.~Eichelsbacher.
\newblock Moderate deviations via cumulants.
\newblock {\em J. Theoret. Probab.}, 26:360--385, 2013.

\bibitem[DJS22]{doering}
H.~D{\"o}ring, S.~Jansen, and K.~Schubert.
\newblock The method of cumulants for the normal approximation.
\newblock {\em Probab. Surv.}, 19:185--270, 2022.

\bibitem[ER59]{ER}
P.~Erd{\Horig{o}}s and A.~R\'enyi.
\newblock On random graphs. {I}.
\newblock {\em Publ. Math. Debrecen}, 6:290--297, 1959.

\bibitem[ER23]{rednos}
P.~Eichelsbacher and B.~Redno{\ss}.
\newblock Kolmogorov bounds for decomposable random variables and subgraph
  counting by the {S}tein-{T}ikhomirov method.
\newblock {\em Bernoulli}, 29(3):1821--1848, 2023.

\bibitem[ET14]{eichelsbacher}
P.~Eichelsbacher and C.~Th{\"a}le.
\newblock New {B}erry-{E}sseen bounds for non-linear functionals of {P}oisson
  random measures.
\newblock {\em Electron. J. Probab.}, 19:no. 102, 25, 2014.

\bibitem[FMN16]{feray}
V.~F\'{e}ray, P.-L. M\'{e}liot, and A.~Nikeghbali.
\newblock {\em Mod-{$\phi$} Convergence}.
\newblock SpringerBriefs in Probability and Mathematical Statistics. Springer,
  Cham, 2016.

\bibitem[GHH95]{hipp}
F.~G{\"o}tze, L.~Heinrich, and C.~Hipp.
\newblock m-dependent random ﬁelds with analytic cumulant generating
  function.
\newblock {\em Scand. J. Statist.}, 22:183--195, 1995.

\bibitem[Gil59]{G}
E.N. Gilbert.
\newblock Random graphs.
\newblock {\em Ann. Math. Statist}, 30(4):1141--1144, 1959.

\bibitem[GT18a]{grotethale18}
J.~Grote and C.~Th{\"a}le.
\newblock Concentration and moderate deviations for {P}oisson polytopes and
  polyhedra.
\newblock {\em Bernoulli}, 24:2811--2841, 2018.

\bibitem[GT18b]{thale18}
J.~Grote and C.~Th{\"a}le.
\newblock Gaussian polytopes: a cumulant-based approach.
\newblock {\em J. Complexity}, 47:1--41, 2018.

\bibitem[GT21]{gusakova}
A.~Gusakova and C.~Th{\"a}le.
\newblock The volume of simplices in high-dimensional {P}oisson–{D}elaunay
  tessellations.
\newblock {\em Ann. H. Lebesgue}, 267:121--153, 2021.

\bibitem[Hei07]{heinrich2}
L.~Heinrich.
\newblock An almost-{M}arkov-type mixing condition and large deviations for
  {B}oolean models in the line.
\newblock {\em Acta Appl. Math.}, 96:247--262, 2007.

\bibitem[HS09]{heinrich}
L.~Heinrich and M.~Spiess.
\newblock Berry-{E}sseen bounds and {C}ram{\'e}r-type large deviations for the
  volume distribution of {P}oisson cylinder processes.
\newblock {\em Lith. Math. J.}, 49:381--398, 2009.

\bibitem[Jan88]{Janson1988}
S.~Janson.
\newblock Normal convergence by higher semiinvariants with applications to sums
  of dependent random variables and random graphs.
\newblock {\em Ann. Probab.}, 16(1):305--312, 1988.

\bibitem[Jan19]{jansen}
S.~Jansen.
\newblock Cluster expansions for {G}ibbs point processes.
\newblock {\em Adv. in Appl. Probab.}, 51(4):1129--1178, 2019.

\bibitem[Kho08]{khorunzhiy}
O.~Khorunzhiy.
\newblock On connected diagrams and cumulants of {E}rd{\Horig{o}}s-{R}\'enyi
  matrix models.
\newblock {\em Comm. Math. Phys.}, 282:209--238, 2008.

\bibitem[KRT17]{reichenbachsAoP}
K.~Krokowski, A.~Reichenbachs, and C.~Th{\"a}le.
\newblock Discrete {M}alliavin-{S}tein method: {B}erry-{E}sseen bounds for
  random graphs and percolation.
\newblock {\em Ann. Probab.}, 45(2):1071--1109, 2017.

\bibitem[LNS21]{LNS21}
G.~Last, F.~Nestmann, and M.~Schulte.
\newblock The random connection model and functions of edge-marked {P}oisson
  processes: second order properties and normal approximation.
\newblock {\em Ann. Appl. Probab.}, 31(1):128--168, 2021.

\bibitem[LP18]{LastPenrose17}
G.~Last and M.~Penrose.
\newblock {\em Lectures on the {P}oisson Process}, volume~7 of {\em Institute
  of Mathematical Statistics Textbooks}.
\newblock Cambridge University Press, Cambridge, 2018.

\bibitem[LPS16]{lastpeccatipenrose}
G.~Last, G.~Peccati, and M.~Schulte.
\newblock Normal approximation on {P}oisson spaces: {M}ehler's formula, second
  order {P}oincar\'e inequality and stabilization.
\newblock {\em Probab. Theory Related Fields}, 165(3-4):667--723, 2016.

\bibitem[LRP13]{lachiezerey2}
R.~Lachi\`eze-Rey and G.~Peccati.
\newblock Fine {G}aussian fluctuations on the {P}oisson space {II}: rescaled
  kernels, marked processes and geometric {$U$}-statistics.
\newblock {\em Stochastic Process. Appl.}, 123(12):4186--4218, 2013.

\bibitem[LRR16]{lachieze-rey}
R.~Lachi\`eze-Rey and M.~Reitzner.
\newblock {$U$}-statistics in stochastic geometry.
\newblock In G.~Peccati and M.~Reitzner, editors, {\em Stochastic Analysis for
  {P}oisson Point Processes: {M}alliavin Calculus, {W}iener-{I}t{\^o} Chaos
  Expansions and Stochastic Geometry}, volume~7 of {\em Bocconi \& Springer
  Series}, pages 229--253. Springer, Berlin, 2016.

\bibitem[LRSY19]{lachiezerey4}
R.~Lachi\`eze-Rey, M.~Schulte, and J.E. Yukich.
\newblock Normal approximation for stabilizing functionals.
\newblock {\em Ann. Appl. Probab.}, 29(2):931--993, 2019.

\bibitem[MM91]{MalyshevMinlos91}
V.A. Malyshev and R.A. Minlos.
\newblock {\em Gibbs Random Fields}, volume~44 of {\em Mathematics and its
  Applications (Soviet Series)}.
\newblock Kluwer Academic Publishers Group, Dordrecht, 1991.

\bibitem[Pen03]{penrosebk}
M.~Penrose.
\newblock {\em Random Geometric Graphs}, volume~5 of {\em Oxford Studies in
  Probability}.
\newblock Oxford University Press, 2003.

\bibitem[Pri12]{momentpoi}
N.~Privault.
\newblock Moments of {P}oisson stochastic integrals with random integrands.
\newblock {\em Probab. Math. Stat.}, 32(2):227--239, 2012.

\bibitem[Pri19]{prkhp}
N.~Privault.
\newblock Moments of $k$-hop counts in the random-connection model.
\newblock {\em J. Appl. Probab.}, 56(4):1106--1121, 2019.

\bibitem[Pri24]{privaultkhops}
N.~Privault.
\newblock Asymptotic analysis of $k$-hop connectivity in the 1{D} unit disk
  random graph model.
\newblock {\em Methodol. Comput. Appl. Probab.}, 26(47), 2024.

\bibitem[PS20]{PS2}
N.~Privault and G.~Serafin.
\newblock Normal approximation for sums of discrete {$U$}-statistics -
  application to {K}olmogorov bounds in random subgraph counting.
\newblock {\em Bernoulli}, 26(1):587--615, 2020.

\bibitem[PS22]{PS4}
N.~Privault and G.~Serafin.
\newblock Berry-{E}sseen bounds for functionals of independent random
  variables.
\newblock {\em Electron. J. Probab.}, 27:1--37, 2022.

\bibitem[PT11]{peccatitaqqu}
G.~Peccati and M.~Taqqu.
\newblock {\em Wiener Chaos: Moments, Cumulants and Diagrams: A survey with
  Computer Implementation}.
\newblock Bocconi \& Springer Series. Springer, 2011.

\bibitem[PY01]{penrose01}
M.~Penrose and J.E. Yukich.
\newblock Central limit theorems for some graphs in computational geometry.
\newblock {\em Ann. Appl. Probab.}, 11(4):1005--1041, 2001.

\bibitem[PY05]{penrose05}
M.~Penrose and J.E. Yukich.
\newblock Normal approximation in geometric probability.
\newblock In {\em Stein's method and applications}, volume~5 of {\em Lect.
  Notes Ser. Inst. Math. Sci. Natl. Univ. Singap.}, pages 37--58. Singapore
  Univ. Press, Singapore, 2005.

\bibitem[R{\"o}l22]{roellin2}
A.~R{\"o}llin.
\newblock Kolmogorov bounds for the normal approximation of the number of
  triangles in the {E}rd{\Horig{o}}s-{R}\'enyi random graph.
\newblock {\em Probab. Engrg. Inform. Sci.}, 36(3):747--773, 2022.

\bibitem[Rot64]{rota1964}
G.-C. Rota.
\newblock On the foundations of combinatorial theory. {I}. {T}heory of
  {M}\"obius functions.
\newblock {\em Z. Wahrscheinlichkeitstheorie und Verw. Gebiete}, 2:340--368,
  1964.

\bibitem[RS13]{reitzner}
M.~Reitzner and M.~Schulte.
\newblock Central limit theorems for ${U}$-statistics of {P}oisson point
  processes.
\newblock {\em Ann. Probab.}, 41(6):3879--3909, 2013.

\bibitem[RSS78]{rudzkis}
R.~Rudzkis, L.~Saulis, and V.A. Statulevi\v{c}ius.
\newblock A general lemma on probabilities of large deviations.
\newblock {\em Litovsk. Mat. Sb.}, 18(2):99--116, 217, 1978.

\bibitem[Ruc88]{rucinski}
A.~Ruci{\'n}ski.
\newblock When are small subgraphs of a random graph normally distributed?
\newblock {\em Probab. Theory Related Fields}, 78:1--10, 1988.

\bibitem[Sch16]{schulte}
M.~Schulte.
\newblock Normal approximation of {P}oisson functionals in {K}olmogorov
  distance.
\newblock {\em J. Theoret. Probab.}, 29:96--117, 2016.

\bibitem[SS91]{saulis}
L.~Saulis and V.A. Statulevi\v{c}ius.
\newblock {\em Limit Theorems for Large Deviations}, volume~73 of {\em
  Mathematics and its Applications (Soviet Series)}.
\newblock Kluwer Academic Publishers Group, Dordrecht, 1991.

\bibitem[ST24]{schulte-thaele}
M.~Schulte and C.~Th{\"a}le.
\newblock Moderate deviations on {P}oisson chaos.
\newblock {\em Electron. J. Probab.}, 29:1--27, 2024.

\bibitem[Zha22]{zhangzs}
Z.S. Zhang.
\newblock Berry-{E}sseen bounds for generalized {$U$}-statistics.
\newblock {\em Electron. J. Probab.}, 27:1--36, 2022.

\end{thebibliography}
\end{document}